%% file: main.tex
\documentclass[11pt, a4paper]{article}
\pdfoutput=1
\usepackage{graphicx}        
                             
\usepackage{algpseudocode}
\usepackage{algorithm}

\usepackage{multicol}        
\usepackage[bottom]{footmisc}

\usepackage{amsmath,amsfonts,amssymb}
\usepackage[varvw]{newtxmath}       
\usepackage[margin=2.5cm]{geometry}
\usepackage{hyperref}

\usepackage{url}
\usepackage{adjustbox}
\usepackage{multirow}

\usepackage{booktabs} 
\usepackage{xcolor}

\usepackage{authblk}

\DeclareMathOperator{\prox}{prox}
\DeclareMathOperator{\shrink}{shrink}

\DeclareMathOperator{\proj}{proj}
\DeclareMathOperator*{\argmin}{arg\,min}

\newcommand{\A}{\mathrm{A}}

\newcommand{\K}{\mathrm{K}}

\newcommand{\mD}{\mathrm{D}}
\newcommand{\mF}{\mathrm{F}}
\newcommand{\TV}{\mathrm{TV}}

\newcommand{\mH}{\mathrm{H}}
\newcommand{\R}{\mathbb{R}}

\usepackage{autonum}
\usepackage{bm}

\begin{document}

\title{Bilevel learning of regularization models and their discretization for image deblurring and super-resolution}

\author[1]{Tatiana A.~Bubba}
\affil[1]{Department of Mathematical Sciences, University of Bath}

\author[2]{Luca Calatroni}
\affil[2]{CNRS, UCA, INRIA, Laboratoire I3S}

\author[3]{Ambra Catozzi$^*$}
\affil[3]{Department of Mathematical, Physical and Computer Sciences, University of Parma}

\author[4]{Serena Crisci$^*$}
\affil[4]{Department of Mathematics and Physics, University of Campania ``Luigi Vanvitelli''}

\author[5]{Thomas Pock}
\affil[5]{Institute of Computer Graphics and Vision, Graz University of Technology}

\author[6]{Monica Pragliola$^*$}
\affil[6]{Department of Mathematics and Applications, University of Naples Federico II}

\author[7]{Siiri Rautio}
\affil[7]{Department of Mathematics and Statistics, University of Helsinki}

\author[8]{Danilo Riccio$^*$}
\affil[8]{School of Mathematical Sciences, Queen Mary University of London}

\author[9]{Andrea Sebastiani$^*$}
\affil[9]{Department of Mathematics, University of Bologna}
\date{}
\maketitle

\abstract{Bilevel learning is a powerful optimization technique that has extensively been employed in recent years to bridge the world of model-driven variational approaches with data-driven methods. Upon suitable parametrization of the desired quantities of interest (e.g., regularization terms or discretization filters), such approach computes optimal parameter values by solving a nested optimization problem where the variational model acts as a constraint. In this work, we consider two different use cases of bilevel learning for the problem of image restoration. First, we focus on learning scalar weights and convolutional filters defining a Field of Experts regularizer to restore natural images degraded by blur and noise. For improving the practical performance, the lower-level problem is solved by means of a gradient descent scheme combined with  a line-search strategy based on the Barzilai-Borwein rule. As a second application, the bilevel setup is employed for learning a discretization of the popular total variation regularizer for solving image restoration problems (in particular, deblurring and super-resolution). Numerical results show the effectiveness of the approach and their generalization to multiple tasks.
}
{\let\thefootnote\relax\footnote{{$^*$ AC, SC, MP, DR and AS contributed equally.}}}

\section{Introduction}\label{sec:intro}

During the last decade, bilevel learning approaches have extensively been used in the field of imaging and vision, see \cite{crockett2022bilevel} for an extensive survey. Originally formulated both in discrete \cite{kunisch2013bilevel} and infinite-dimensional \cite{de2013image} settings as \emph{shallow} learning strategies for estimating the parametrized image regularizers \cite{ochs2015bilevel} and noise models \cite{calatroni2017bilevel,Calatroni_2019}, they have further been extended to more challenging scenarios involving higher-order regularization models \cite{TuomoJMIV2017}, space-adaptive regularization \cite{DeLosReyes2022}, non-local extensions \cite{DElia2021} and optimal discretization \cite{chambolle2021learning}. Deep variants of bilevel models can naturally be considered by means of algorithmic unrolling of iterative solvers \cite{Bonettini2022,Lorenz2022} and interesting connections with Deep Equilibrium Models can also be shown \cite{franceschi2018bilevel,Riccio2022}. From an optimization viewpoint, bilevel learning  is becoming an increasingly popular approach for interpreting automatic differentation techniques, see, e.g., \cite{ji2021bilevel,ehrhardt2021inexact,Bolte2022}.

In this work, we consider bilevel learning approaches for estimating tailored regularization models and their discretization in the framework of ill-posed linear inverse problems \cite{engl1996regularization,hansen2006deblurring} of the form
\begin{equation}
  \label{eq:lin_mod}  \text{find}\qquad\bm{u}\in\R^n\qquad\text{such that}\qquad \R^m\ni \bm{f}=\bm{\A}\bm{u}+\bm{e}\,,
\end{equation}
where  $m\leq n$ and $\bm{u}$ is the vectorized unknown image to retrieve from noisy, blurred and possibly under-sampled data $\bm{f}$, $\bm{e}\in\R^m$ is an additive (e.g., white Gaussian) noise component and $\bm{\A}\in\R^{m\times n}$ is a known forward model. Notable examples considered in this work will be the two cases $\bm{\A} = \mathbf{H}\in\R^{n\times n}$, a blur matrix, and $\bm{\A} = \mathbf{SH} \in\R^{m\times n}$ where $\mathbf{S}\in\R^{m\times n}$ is a decimation operator.

To overcome the ill-posedness often arising for problems \eqref{eq:lin_mod}, a variational framework is often considered (see, e.g., \cite{benning2018modern} for a survey) where an approximate solution $\bm{u}^*\in\R^{n}$ is defined as the minimizer of a parametrized functional $\mathcal{J}$ combining prior knowledge on the acquisition process modeled in a data fidelity term with some \emph{a priori} information on the unknown $\bm{u}$.  
While the choice of the data fidelity term is often hand-crafted by following suitable statistical considerations based on Maximum A Posteriori estimation \cite{bayes}, a customized choice of a suitable image regularizer to the particular application and setup considered is often challenging.
Taking as an example the Total Variation (TV) regularization \cite{rof}, it is nowadays common knowledge that for a better and more flexible modeling of the quantities of interest, hyper-parametrized regularization models should be used \cite{holler2018bilevel}.

As a result, in the last years a lot of attention has been devoted to the design of efficient, robust and automatic strategies for the parameters identification task. In this perspective, besides statistical \cite{Vidal2020} and deep learning \cite{Afkham_2021} techniques, bilevel approaches tackle the task as a nested optimization problem formulated in terms of a loss functional $\Phi$ defined in terms of a training dataset composed of pairs $\{(\bm{g}_j,\bm{f}_j)\}_{j=1}^s$ with $\bm{g}_j\in\R^n$ and $\bm{f}_j\in\R^m$ denoting ground truth images and their corresponding degraded acquisitions, respectively. {In formulas, the problem reads:}
\begin{align}\label{eq:bil_or}
 &\min_{\bm{\theta}}~\frac{1}{s}\sum_{j=1}^s \Phi(\bm{u}_j(\bm{\theta}),\bm{g}_j)   \\
 &\text{subject to}\quad \bm{u}_j(\bm{\theta})\in\argmin_{\bm{u}}~\mathcal{J}(\bm{u};\bm{\theta},\bm{f}_j)\quad j=1,\ldots,s
\end{align}
where $\bm{\theta}\in\R^{d}, d\gg 1$ codifies the unknown parameters to estimate,  $\Phi:\R^n\times \R^n\to\R_{+}$ is a loss function and $\mathcal{J}:\R^n\times\R^m\times\R^d\to \R_+$ is the reconstruction model acting here as a constraint. 

In this work, we consider two different case studies under specific assumptions for problem \eqref{eq:bil_or}. In both cases we assume that the noise arising in the degradation model \eqref{eq:lin_mod} is additive white Gaussian so that, for a given $\bm{f}\in\R^m$ the cost functional $\mathcal{J}$ can be specified as
\begin{equation}\label{eq:min_functional}
    \mathcal{J}(\bm{u};\bm{\theta},\bm{f}) = \mathcal{R}(\bm{u};\bm{\theta}) + \frac{1}{2}\|\bm{\A u}(\bm{\theta})-\bm{f}\|_2^2\, 
\end{equation}
and where $\mathcal{R}(\bm{u};\bm{\theta})$ will be specified in two different instances.
Namely, in the first case study we consider the estimation of the parameters $\bm{\theta}=(\alpha_1,\ldots,\alpha_L, \bm{k_1},\ldots,\bm{k}_L)\in (\R_+)^L\times (\R^{\kappa\times \kappa})^L$ of a Field of Experts (FoE) regularization term \cite{foe,kunisch2013bilevel,pock1} defined by
\begin{equation}  \label{eq:foe}
    \mathcal{R}(\bm{u};\bm{\theta}) = \sum_{\ell=1}^L \sum_{i=1}^n \alpha_\ell\varphi((\bm{k}_{\ell}\ast \bm{u})_i)\,,\quad \varphi(x) =  \log(1+|x|^2)\,,
\end{equation}
where $\ast$ denotes the convolution product, $\alpha_1,\ldots,\alpha_L$ are positive weights and $\bm{k}_1,\ldots\bm{k}_L$ are convolution kernels of size $\kappa\times \kappa$. 
In \cite{pock1} a bilevel problem in the form \eqref{eq:bil_or} with \eqref{eq:min_functional}-\eqref{eq:foe} is studied with $\bm{\A}=\mathbf{I}$, i.e. for image denoising problems. The estimated regularized model is then tested on more complex imaging problems (such as image deblurring) showing good generalization properties. In the following, we include more explicitly some degradation models $\bm{\A}\neq \mathbf{I}$ in order to better adapt the regularization model to actual image restoration problems.
As an efficient lower-solver improving the overall efficiency of the bilevel scheme, we employ a line-search gradient method combined with Barzilai-Borwein (BB) steplength updating strategies~\cite{barzilai1988two}.
This choice is, indeed, particularly appealing since BB rules can be easily extended to gradient schemes for efficiently solving general (non-convex) optimization problems, still preserving their low memory requirements and low computational cost per iteration \cite{raydan1997barzilai,di2018steplength}. 

\medskip

In the second case study, we focus on the optimal selection via bilevel learning of discretization filters for the popular TV regularization term. This is indeed a crucial question in practical applications. As shown in \cite{condat2017discrete,Hintermuller_2014}, in order to prevent discretization biases in the reconstructed images, suitable finite difference discretization stencils (expressed here as convolution filters) may be designed. In \cite{chambolle2021learning} an analogous strategy was employed: interestingly, it was shown that the optimal filters learned on a specific task (there, denoising/inpainting) do not always generalize well when applied to other tasks. Based on this observation, we proposed in the following an analogous bilevel strategy for learning TV discretization filters for more challenging problems: deblurring and super-resolution.
For doing so and denoting by $\bm{\mD}:\mathbb{R}^{\sqrt{n}\times\sqrt{n}}\rightarrow\mathbb{R}^{n}\times\mathbb{R}^{n}$ , $\bm{\mD}\bm{u}=(\bm{\mD}^v\bm{u}, \bm{\mD}^h\bm{u}) \in \R^{n\times n}$ the discrete gradient operator acting on the 2D image\footnote{We consider square images of size $\sqrt{n}\times\sqrt{n}$ for simplicity, with $\sqrt{n}\in\mathbb{N}$.} $\bm{u}$,
we introduce a parameterized discrete version of TV defined in terms of an averaging operator $\bm{\mF}:\mathbb{R}^{n}\times\mathbb{R}^{n}\rightarrow\bm{Z}$  with $\bm{Z}$ being Cartesian product of $L\ge 1$ copies of $\mathbb{R}^{n}\times\mathbb{R}^{n}$ with norm
\begin{equation}  \label{norm:dual}
    \|(\bm{z}^1,\ldots,\bm{z}^L)\|_{\bm{Z}}:=\sum_{i=1}^{L}\|\bm{z}^i\|_{1,2}.
\end{equation}
The generalized dual formulation of TV (see \cite{Hintermuller_2014,condat2017discrete}) is thus defined by:
\begin{equation} \label{eq:dual_TV}
    \mathcal{R}(\bm{u};\bm{\mF}) = \TV_{\bm{\mF}}(\bm{u}):=\sup\{\langle\bm{p},\bm{\mD}\bm{u}\rangle:\ \|\bm{\mF}\bm{p}\|_{\bm{Z}}^*\le 1\}
\end{equation}
where $\bm{p}=(\bm{p}^1,\bm{p}^2)\in\R^{n\times n}$ denotes the discrete dual variables, $\|\cdot\|^*_Z$ denotes the dual norm defined on $Z^*$. For each $l=1,\ldots, L$, the two components of the operator $\bm{\mF}^l= (\bm{\mF}^{l,1},\bm{\mF}^{l,2})$ acts as convolutions on $\bm{p}^1$ and $\bm{p}^2$, respectively, thus averaging appropriately their discretized values on the grid.
Note, in particular, that both \eqref{norm:dual} and \eqref{eq:dual_TV} are indeed generalized version of the usual definition of TV which in its standard form corresponds to the case $L=1$, $\bm{z}=\bm{\mD u}$ and $(\bm{\mF p})_{i,j} = (p^1_{i+\frac{1}{2},j}, p^2_{i,j+\frac{1}{2}})$, with:
\[
\| \bm{\mF p} \|^*_Z = \sum_{i,j=1}^n \sqrt{ \left(p^1_{i+\frac{1}{2},j}\right)^2 + \left(p^2_{i,j+\frac{1}{2}}\right)^2},
\]
upon identification of $\bm{Z}^*$ with $\bm{Z}$.
The problem of learning the optimal TV discretization filters takes here the form:
\begin{align}\label{eq:bil_2}
 &\min_{\bm{\mF}}~\frac{1}{s}\sum_{j=1}^s \Phi(\bm{u}_j(\bm{\mF}),\bm{g}_j)   \\
 &\text{s.t.}\ \bm{u}^*_j(\bm{\mF})\in\argmin_{\bm{u},\bm{q}}\max_{\bm{p}}\langle\bm{\mD}\bm{u}-\bm{\mF}^*\bm{q}, \bm{p}\rangle + \lambda\|\bm{q}\|_{Z}+\frac{1}{2}\|\bm{\A u}-\bm{f}_j\|_2^2\quad j=1,\ldots,s,
\end{align}
where a primal-dual formulation is used as lower-level constraint.
The solution of both the lower-level and the nested bilevel problem requires here some attention. As far as the lower-level problem is concerned, given its primal-dual structure, a natural choice for computing approximate solutions would be considering a (preconditioned) primal-dual algorithm \cite{ChambollePock2011,PockChambolle2011_prec}. When embedded in a bilevel learning framework, however, such strategy may be not ideal. In order to compute derivatives of the upper-level problem with respect to $\bm{\mF}$ by means of automatic differentiation an back-propagation, one could indeed run a primal-dual solver for a sufficiently high number of inner iterations to guarantee a good numerical approximation. Depending on the problem at hand, however, this number of iterations could be very high (the harder the problem, the higher the number, heuristically) which could cause memory shortages. Furthermore, to compute gradient updates of the upper-level problem, one should also compute the adjoint states associated with the bilevel problem, which adds computational burden to the overall optimization. Following recent works \cite{Bogensperger2022,chambolle2021learning}, we will consider in this work a piggy-back primal-dual algorithm computing solutions of the lower-level problem and of the adjoint states at the same time. We consider both the problem of image deblurring ($\bm{\A}=\mathbf{H}\in\R^{n\times n}$, a structured circulant convolution matrix) and super-resolution ($\bm{\A}=\mathbf{SH}\in\R^{m\times n}$) and, in the latter case, we resort to Fourier-based approaches previously proposed in \cite{zhao2016fast} and used, e.g., in \cite{Pragliola2023, mylonopoulos2022constrained} for computing proximal updates in a closed-form.

\medskip

We remark that while the two case studies \eqref{eq:bil_or} with \eqref{eq:min_functional}-\eqref{eq:foe} and \eqref{eq:bil_2} share several analogies due to the common bilevel optimization framework considered, they are intrinsically different from an optimization view point. While both bilevel problems are in fact globally non-convex, in the former case the lower-level problem is non-convex and smooth, while in the latter case the lower-level problem is convex and non-smooth. It is therefore natural to exploit smoothness in the former case by designing an efficient gradient-type solver (relying on BB line-search) and convexity in the latter case by means of appropriate primal-dual updates.

\medskip

\noindent \emph{Structure of the paper.} In Section \ref{sec:foe} we detail the bilevel approach for estimating optimal parameters $\bm{\theta}$ of a FoE regularization model \eqref{eq:foe} for image restoration problems. Similarly, in Section \ref{sec:TV_disc} we detail the bilevel strategy for learning optimal total variation discretization filters in image restoration tasks. Several numerical results on both tasks are reported in Section \ref{sec:experiments}. Closing remarks are given in Section~\ref{sec:conclusions}.
\input{part_foe}
\input{part_tv}

\section{Numerical experiments}\label{sec:experiments}
\input{part_results}

\section{Conclusions}\label{sec:conclusions}
In this work, we applied the bilevel learning framework presented in  \cite{pock1,chambolle2021learning} for learning image regularization models and their discretization to some exemplar image restoration problems (deblurring, super-resolution). 

Starting from~\cite{pock1}, we adapted the regularization model to explicitly account for the blur operator, by including it in the formulation of the lower-level problem. In addition, to exploit  the smoothness  of the lower-level problem, we equipped the gradient descent scheme with a non-monotone line-search coupled with the BB rule. Numerical results show that even with a small training set of only 9 images, trained on different types of blur, the proposed approach delivers reliable results: noise is suppressed and, with the exception of motion blur which remains a more challenging scenario, blur distortions are effectively removed. 

In the second part, we proposed a strategy analogous to~\cite{chambolle2021learning} for learning TV discretization filters for two more challenging image restoration problems: deblurring and super-resolution, and for both we trained with noisefree and noisy data. This required to update the computation of the proximal operator in the piggy-back scheme, to account for the presence of a non-trivial forward operator. Numerical results are consistent with those in~\cite{chambolle2021learning}.

To improve the generalization capabilities, it would be worth augmenting the number and type of degradations, by considering the problem in a multi-task learning framework.

\section*{Acknowledgements}
The authors would like to thank the organizers and participants of the workshop on Advanced Techniques in Optimization for Machine learning and Imaging (ATOMI, Rome, 20-24 June, 2022) during which the present work was initiated. This work has been partially supported by the INDAM-GNCS research group.
This research utilized Queen Mary's Apocrita and Andrena HPC facilities, supported by QMUL Research-IT \url{http://doi.org/10.5281/zenodo.438045}.
This research also utilized the HPC facilities from the Department of Mathematics of the University of Bologna.
S. Crisci and M. Pragliola aknowledge the support from the EU-FESR PON Ricerca e Innovazione 2014-2020, art. 24, comma 3, lett. a) L. 240/2010 e s.m.i., D.M. 1062/2021.

\bibliography{references}{}
\bibliographystyle{plain}

\end{document}

%% file: part_foe.tex
\section{Bilevel learning of FoE regularization}
\label{sec:foe}

Given a training set composed by $s$ couples of images $\{(\bm{g}_j,\bm{f}_j)\}_{j=1}^s$, 
we consider an instance of \eqref{eq:bil_or} for the estimation of the weights $\bm{\alpha}=(\alpha_1,\ldots,\alpha_L)$ and the convolution filters $\bm{k}=(\bm{k}_1,\ldots,\bm{k}_L)$ defining the Field of Expert (FoE) defined in \eqref{eq:min_functional}-\eqref{eq:foe}. For simplicity, we consider a quadratic loss function for the upper level problem. The bilevel problem reads:
\begin{align}\label{eq:bil_foe}
 &\min_{\bm{\theta}\in\mathcal{C}}~ \left\{ \frac{1}{s}\sum_{j=1}^s \Phi(\bm{u}_j(\bm{\theta}),\bm{g}_j) =  ~\frac{1}{2s}\sum_{j=1}^s \|\bm{u_j}(\bm{\theta})-\bm{g_j}\|_2^2 \right\} \\
 &\text{s.t.}\quad \bm{u}_j(\bm{\theta})\in\argmin_{\bm{u}}~\left\{ \mathcal{J}(\bm{u};\bm{\theta},\bm{f}_j) = \frac{1}{2}\|\bm{\mH}\bm{u}-\bm{f_j}\|_2^2 +\sum_{\ell=1}^L\sum_{i=1}^n\alpha_\ell\varphi((\bm{k}_\ell\ast \bm{u})_i)\right\}
\end{align}
for all $j=1,\ldots,s$, where, for all $\ell=1,\ldots,L$ we recall
\begin{equation}
    \label{eq:phi}
    \varphi((\bm{k}\ast\bm{u})_i) =  \log(1+|(\bm{k}_{\ell}\ast \bm{u})_i|^2)\quad i=1,\ldots, n
\end{equation}
and where $\mathcal{C}$ is the constraint set for the unknowns $\bm{\theta}=(\bm{\alpha},\bm{k})$ defined by
\begin{equation}\label{eq:con_set}
    \mathcal{C}=C_{\bm{\alpha}}\times C_{\bm{k}}\,,\;C_{\bm{\alpha}}=\R_+^L\,,\;C_{\bm{k}} = \varprod_{\ell=1}^L C_{\bm{k}_{\ell}}\,,\;C_{\bm{k}_{\ell}}:=\{\bm{k}_\ell\mid\bm{1}^T\bm{k}_\ell=0\}\,,
\end{equation}
where $\R_+$ denotes the set of non-negative real numbers.
To compute the gradients of the loss function with respect to $\bm\theta$, we exploit the smoothness of the lower-level optimization problems by replacing them with their first-order optimality condition. By expressing the convolution products as matrix-vector products, so that $\bm{k}_\ell\ast\bm{u}=\bm{\K}_\ell\bm{u}$, we thus consider consider the following problem:
\begin{align}
\label{eq:bil_foe2}
\begin{split}
&\min_{\bm{\theta}\in \mathcal{C}} ~
\frac{1}{2s}\sum_{j=1}^s \|\bm{u_j}(\bm{\theta})-\bm{g_j}\|_2^2 \\
&\text{s.t. }\nabla_{\bm{u}}\mathcal{J}(\bm{u}_j(\bm{\theta}),\bm{\theta})=\bm{0},\quad  j=1,\ldots,s,
\end{split}
\end{align}
where, by dropping the dependence on $\bm{\theta}$ for all $\bm{u}_j=\bm{u}_j(\bm{\theta})$  for ease of notation, the optimality conditions read as:
\begin{equation} \label{eqn:opt_cond_gradJ}
  \nabla_{\bm{u}}\mathcal{J}(\bm{u}_j,\bm{\theta}) = \bm{\mH}^T(\bm{\mH}\bm{u}_j-\bm{f}_j)+\sum_{\ell=1}^{L}\sum_{i=1}^n\alpha_\ell\bm{\K}_\ell^T\varphi'((\bm{\K}_\ell\bm{u}_j)_i) = \bm 0,\quad j=1,\ldots,s.
\end{equation}
The Lagrangian functional associated with the above constrained minimization problem reads:
\begin{align}
     \mathcal{L}(\bm{u},\bm{\alpha},\bm{k},\bm{p},\bm{q},\bm{r})  = &\frac{1}{2s}\sum_{j=1}^s \|\bm{u_j}-\bm{g_j}\|_2^2\;{+}\sum_{j=1}^s\;\left\langle\sum_{\ell=1}^{L}\sum_{i=1}^n\alpha_\ell\bm{\K}_\ell^T\varphi'((\bm{\K}_\ell\bm{u}_j)_i),\bm{p}_j \right\rangle  \\&
      \;{+}\;\sum_{j=1}^s\left\langle\bm{\mH}^T(\bm{\mH}\bm{u}_j-\bm{f}_j),\bm{p}_j\right\rangle
 \;{-}\; \sum_{\ell=1}^{L} 
 \left( q_\ell(\bm{1}^T\bm{k}_\ell) + \alpha_\ell r_\ell \right)\,,
\end{align}
where $\bm{p}_j\in\R^n$ and $\bm{q},\bm{r}\in\R^L$ are the vectors of the Lagrange multipliers associated with the constraint in \eqref{eq:bil_foe2} and in \eqref{eq:con_set}, respectively. 
The first-order necessary optimality conditions for problem~\eqref{eq:bil_foe2} guarantee the existence of the vectors
 $\bm{p}_1,\ldots,\bm{p}_s,\bm{q},\bm{r}$ such that for $j=1,\ldots,s$ we have:

\begin{align}\label{eq:opt1}
&\left(\sum_{\ell=1}^L\sum_{i=1}^n\alpha_\ell\bm{\K}_\ell^T \bm{\K}_\ell\varphi''((\bm{\K}_\ell\bm{u}_j)_i) + \bm{\mH}^T\bm{\mH}\right)\bm{p}_j+\bm{u}_j-\bm{g}_j\;{=}\;0  \\
     \label{eq:opt2}
&\Bigg\langle\sum_{i=1}^n\bm{\K}_\ell^T\varphi'((\bm{\K}_\ell\bm{u}_j)_i),\bm{p}_j\Bigg\rangle-r_\ell=0\,, \quad \ell=1,\ldots,L
 \\
     \label{eq:opt3}
& \nabla_{\bm{\K}_\ell}\Bigg\langle\sum_{i=1}^n\alpha_\ell\bm{\K}_\ell^T\varphi'((\bm{\K}_\ell\bm{u}_j)_i),\bm{p}_j\Bigg\rangle-q_\ell =0\,,\quad \ell=1,\ldots,L\\
    \label{eq:opt4}
    &\bm{\mH}^T(\bm{\mH}\bm{u}_j-\bm{f}_j)+\sum_{\ell=1}^{L}\sum_{i=1}^n\alpha_{\ell}\bm{\K}_\ell^T\varphi'((\bm{\K}_\ell\bm{u}_j)_i) \;{=}\;0  \\ 
  \label{eq:opt5}
  &\bm{1}^T\bm{k}_\ell=0\,,\quad \ell=1,\ldots,L\,\\
  \label{eq:opt6}
  &\bm{r}-\max(\bm{0},\bm{r}-c\bm{\alpha})=0\,.
\end{align}
Equation \eqref{eq:opt6} is a slackness condition equivalent to $\bm{\alpha}\geq 0$, $\bm{r}\geq 0$ and $\langle\bm{\alpha},\bm{r}\rangle=0$ where the $\max$ is meant component-wise and $c\in\R$ is any positive scalar.  Note that the Lagrange multipliers $\bm{q},\bm{r}$ associated with the constraints expressed in \eqref{eq:con_set} do not need to be explicitly computed, as the required properties on the weights~$\bm{\alpha}$ and on the filters~$\bm{k}$ can be directly handled by suitable projections imposed during the optimization used for the computation of~$\bm{\theta}$. We refer the reader to \cite{pock1} for more details.
In our setting, the unknown of primary interest is represented by $\bm{\theta}$. To compute it, we first solve problems \eqref{eq:opt1} and \eqref{eq:opt4} in order to recover the expressions of $\bm{u}_j,\bm{p}_j$ to be plugged into \eqref{eq:opt2} and \eqref{eq:opt3}. 
Problems \eqref{eq:opt2} and \eqref{eq:opt3} provide the gradient components of the upper-level functional with respect to $\bm \theta$, which are then employed within an iterative gradient method applied endowed with a suitable projection onto the constraint set $\mathcal{C}$ for computing the solution. 
Denoting by $ \nabla_{\bm{u}}^2\mathcal{J}(\bm{u}_j,\bm{\theta})$ the Hessian matrix of $\mathcal{J}(\bm{u}_j,\bm{\theta})$ so that:
\[
 \nabla_{\bm{u}}^2\mathcal{J}(\bm{u}_j,\bm{\theta}) = \sum_{\ell=1}^L\sum_{i=1}^n\alpha_\ell\bm{\K}_\ell^T \bm{\K}_\ell\varphi''((\bm{\K}_\ell\bm{u}_j)_i) + \bm{\mH}^T\bm{\mH},\quad j=1,\ldots,s,
\]
we can compute by suitable manipulations of \eqref{eq:opt1}-\eqref{eq:opt4}
\footnotesize
\begin{eqnarray} \label{eqn:grad_phi_alpha}
\left(\nabla_{\bm \alpha}\Phi(\bm{u}_j, \bm{g}_j )\right)_{\ell} &=& \left(\sum_{i=1}^n\bm{\K}_\ell^T\varphi'((\bm{\K}_\ell\bm{u}_j)_i) \right)^T\left(\nabla_{\bm{u}}^2 \mathcal{J}(\bm{u_j},\bm{\theta})\right)^{-1}(\bm g_j - \bm u_j), \\
\label{eqn:grad_phi_kl}
 \left(\nabla_{\bm{k}}\Phi(\bm{u}_j, \bm{g}_j ) \right)_{\ell} &=& \left(\sum_{i=1}^n\alpha_\ell\varphi'((\bm{\K}_\ell\bm{u}_j)_i) + \sum_{i=1}^n\alpha_\ell \bm \K_{\ell}^T \bm \K_{\ell}\varphi''((\bm{\K}_\ell\bm{u}_j)_i)\right)^T \left(\nabla_{\bm{u}}^2  \mathcal{J}(\bm{u_j},\bm{\theta})\right)^{-1}(\bm g_j - \bm u_j), \notag
\end{eqnarray}
\normalsize
for $\ell = 1,\dots, L$. Note that seeking for stationary points of the functional $\mathcal{J}(\bm{u},\bm{\theta})$ with $\bm{\theta}$ fixed can be costly. Alternatively, one can directly address the lower-level minimization problem expressed by \eqref{eq:bil_foe} by computing high-precision solutions at each outer iteration.
Upon a suitable initialization $(\bm{u}^{(0)},\bm{\theta}^{(0)})$, the $k$-th iteration of the scheme for the update of $\bm{\theta}=(\bm{\alpha},\bm{k})$ reads, for all $j=1,\ldots,s$:
\begin{align} 
\label{eq:sub_u}
\bm{u}_j^{(k+1)}\;{\in}\;&\argmin_{\bm{u}}\left\{\frac{1}{2}\|\bm{\mathrm{H} u}-\bm{f}_j\|_2^2+\sum_{\ell=1}^{L}\sum_{i=1}^n\alpha_\ell^{(k)}\varphi((\bm{\K}_\ell^{(k)}\bm{u_j})_i)\right\}\\
\label{eq:sub_p}
\bm{p}_j^{(k+1)} \;{=}\;&\left(\sum_{\ell=1}^L \sum_{i=1}^n\alpha_\ell^{(k)}(\bm{\K}_\ell^{(k)})^T \bm{\K}_\ell^{(k)})\varphi''((\bm{\K}^{(k)}_\ell\bm{u}_j^{(k+1)})_i) + \bm{\mH}^T\bm{\mH}\right)^{-1}(\bm{g}_j-\bm{u}_j^{(k+1)})\\
\label{eq:sub_alpha}
    \bm{\alpha}^{(k+1)} \;{=}\;&P_{C_{\bm{\alpha}}}\left(\bm{\alpha}^{(k)}-\tau \frac{1}{s}
    \nabla_{\bm{\alpha}}\Phi(\bm{u}_j^{(k+1)},\bm{g}_j)\right)\\
        \label{eq:sub_k}
\bm{k}^{(k+1)} \;{=}\;&P_{C_{\bm{k}}}\left(\bm{k}^{(k)}-\tau\frac{1}{s}
\nabla_{\bm{k}}\Phi(\bm{u}_j^{(k+1)},\bm{g}_j)\right)
\end{align}
with $\nabla_{\bm{\alpha}}\Phi, \nabla_{\bm{k}}\Phi$ defined component-wise by equations \eqref{eqn:grad_phi_alpha}. 
Also, notice that the updates of $\bm{\alpha}$, $\bm{k}$ amount to one step of projected gradient descent with step-size $\tau$, with the projection being performed on the sets $C_{\bm{\alpha}}$, $C_{\bm{k}}$ introduced in \eqref{eq:con_set}.
More specifically, the non-negativity constraint for $\alpha_\ell$, $\ell=1,\ldots,L$ can be easily addressed by projecting the updated weights onto $[0,+\infty)$; for what concerns the filters, we subtract to each $\bm{k}_\ell^{(k+1)}$, $\ell=1,\ldots,L$, its mean, so that the condition in \eqref{eq:con_set} is satisfied. Finally, notice that the projection $P_{C_{\bm{\alpha}}},P_{C_{\bm{k}}}$ account for the Lagrange multipliers $\bm{q},\bm{r}$ that, as discussed above, are not explicitly included in the optimization.

To computing $\bm{u}^{(k+1)}$ at each outer iteration $k\geq 0$  in \eqref{eq:sub_u} we employ a gradient method. By dropping the dependence on $j=1,\ldots,s$ and upon a warm-start initialization $\bm{u}_{(0)} = \bm{u}^{(k)}$, the inner iteration loop reads:
\[
\bm u_{(t+1)} = \bm u_{(t)} + \nu_t \bm d_{(t)}, 
\]
where $\bm d_{(t)} = -\nabla_{\bm{u}}\mathcal{J}(\bm{u}_{(t)},\bm{\theta}^{(k)})$, $\nu_t = \beta^{\xi} \gamma_t$, $\beta \in (0,1)$ and $\xi$ is the first non-negative integer such that the following Armijo decrease condition is satisfied:
\begin{equation}\label{eq:armijo}
f(\bm u_{(t)} +  \beta^{\xi} \gamma_t\bm d_{(t)}) \leq f(\bm u_{(t)}) + \sigma  \beta^{\xi} \gamma_t {\nabla_{\bm{u}}\mathcal{J}(\bm{u}_{(t)},\bm{\theta}^{(k)})}^T \bm {d}_{(t)}.
\end{equation}
Here, the step-length $\gamma_t$ is computed in accordance with the spectral rule $\rm{BB}1$ \cite{barzilai1988two}, defined by
\begin{equation}\label{eq:bb1}
\gamma_t = \frac{ \lVert \bm \rho_{(k-1)} \rVert_2^2 }{ \langle {\bm \rho_{(t-1)}}, \bm y_{(t-1)} \rangle }
\end{equation}
where $\bm \rho_{(t-1)} = \bm u_{(t)} - \bm u_{(t-1)}$ and $\bm y_{(t-1)} = \nabla_{\bm{u}}\mathcal{J}(\bm{u}_{(t)},\bm{\theta}^{(k)}) - \nabla_{\bm{u}}\mathcal{J}(\bm{u}_{(t-1)},\bm{\theta}^{(k)})$; this choice is combined with the safeguarding condition $\gamma_t \in [\gamma_{\rm{min}}, \gamma_{\rm{max}}]$, where $0<\gamma_{\rm{min}}< \gamma_{\rm{max}}$. The main steps of the resulting scheme for solving the lower-level problem are outlined in Algorithm \ref{alg:solve_lower}; for the convergence properties we refer to~\cite{raydan1997barzilai}. A pseudo-code of the lower-level solver is reported in Algorithm \ref{alg:solve_lower}.

Solving \eqref{eq:sub_p} amounts to solve $s$ linear system in the variables $\bm{p}_j, j=1,\ldots,s$ which can be done by means of a Krylov-subspace method, such as, e.g., the Conjugate Gradient (CG) method. The CG iterations are stopped as soon as the residual norm is below a selected tolerance. 
Notice that the invertibility of the coefficient matrix in \eqref{eq:sub_p} is not guaranteed \emph{a priori} as it strongly depends on the updated filters. Nonetheless, in the computed examples reported in Section \ref{sec:experiments} we did not observe ill-conditioning; as a safeguard, the coefficient matrix in \eqref{eq:sub_p} may be slightly changed by adding $\varepsilon \bm{\mathrm{I}}$, with $\varepsilon>0$.
The updates \eqref{eq:sub_alpha}, \eqref{eq:sub_k} can be run up till convergence to compute an approximation of the desired FoE parameters $\bm{\theta}=(\bm{\alpha},\bm{k})$ which can then be used for solving test image restoration problems by means of the optimal FoE regularizer computing via the subroutine \texttt{SolveLower}. The stopping criterion 
for the gradient method employed in the $\bm{u}$-update in \eqref{eq:sub_u} is on the relative change between two consecutive collections of restored images
\begin{equation}\label{eq:rel_chg}
 \frac{1}{s}\sum_{j=1}^s\frac{\|\bm{u}_{j,(t+1)}-\bm{u}_{j,(t)}\|_2}{\|\bm{u}_{j,(t)}\|_2}  <\tau_{\mathrm{inner}}\,,
\end{equation}
whereas the outer scheme is stopped when a maximum number of iterations is reached. 
More details on the selection of the parameters and the inner 
tolerance are given in Section~\ref{sec:experiments}.
The overall bilevel learning procedure is outlined in Algorithm \ref{alg:learn_foe}. 

\begin{algorithm}[h]
\caption{Function \texttt{SolveLower}}\label{alg:solve_lower}
\begin{algorithmic}[0]
\State \textbf{Input}:  $\bm{f},\bm{\mathrm{H}},\bm{\theta}^{(k)}, \bm{u}^{(k)}, \beta,\sigma\in (0,1), 0<\gamma_{\rm min}<\gamma_{\rm max}, \gamma_0>0$, $t=0$
\State $\bm{\cdot}$ $\bm u^{(0)} =  \bm u^{(k)}$ \quad {(\emph{warm start})} 
\While {not converging} 
    \State $\bm{\cdot}$ compute $\bm d_{(t)}= -\gamma_t \nabla_{\bm{u}}\mathcal{J}(\bm{u}_{(t)},\bm{\theta}^{(k)})$ 
    \State $\bm{\cdot}$ set $\nu_t = 1$
    \State $\bm{\cdot}$ perform Armijo backtracking procedure  \eqref{eq:armijo} to update $\nu_t$
    \State  $\bm{\cdot}$ compute $\bm u_{(t+1)} = \bm u_{(t)} + \nu_t \bm d_{(t)}$
    \State $\bm{\cdot}$ update $\gamma_{t+1}\in [\gamma_{\rm min}, \gamma_{\rm max}]$ according to \eqref{eq:bb1} 
    \State $\bm{\cdot}$ $t=t+1$
\EndWhile
\State \textbf{Output}: $\bm u^{(k+1)}=\bm u_{(t+1)}$
\end{algorithmic}
\end{algorithm}

\begin{algorithm}[h]
\caption{FoE-bilevel learning for image restoration problem}\label{alg:learn_foe}
\begin{algorithmic}[0]
\State \textbf{Initialization}: Set $\bm{\theta}^{(0)}=(\alpha^{(0)},\bm{k}^{(0)})$, $k=0$, $\tau >0$
\State \textbf{Input}: $\bm{\mathrm{H}}$ and training examples $\{(\bm{g}_j,\bm{f}_j)\}_{j=1}^s$
\While {not converging}
    \State $\bm{\cdot}$ compute $\bm{u}_j^{(k+1)}=$ \texttt{SolveLower}$(\bm{f},\bm{\mathrm{H}},\bm{\theta}^{(k)}, \bm u_j^{(k)})$, $j=1,\ldots, s$
    \State $\bm{\cdot}$ compute $\bm{p}_j^{(k+1)}$ by solving \eqref{eq:sub_p} with CG, $j=1,\ldots, s$
    \State $\bm{\cdot}$ update $\bm{\alpha}^{(k+1)}$ by means of \eqref{eq:sub_alpha}
      \State $\bm{\cdot}$ update $\bm{k}^{(k+1)}$ by means of \eqref{eq:sub_k}
    \State $\bm{\cdot}$ $k=k+1$
\EndWhile
\State \textbf{Output}: optimal parameters  $\bm{\theta}^*=(\bm{\alpha}^{(k+1)},\bm{k}^{(k+1)})$ 
\end{algorithmic}
\end{algorithm}

%% file: part_tv.tex
\section{Bilevel learning of TV discretization}  \label{sec:TV_disc}
We now consider a different use case: the problem of learning a suitable discretization for the Total Variation (TV) regularization. Following \cite{chambolle2021learning}, we formulate the problem as a bilevel learning problem of a suitable loss function defined in terms of a training set$\{(\bm{g}_j,\bm{f}_j)\}_{j=1}^s$ as above. Recalling the dual definition of TV provided in \eqref{eq:dual_TV}, we consider a family of convolution-type discretization operators $\bm{\mF}=(\bm{\mF}^l)_{l=1}^L$ acting on the TV dual variable  for $\ell=1,\ldots, L$ as
    \begin{equation}
        (\bm{\mF}^l\bm{p})=
        \begin{pmatrix}
        \bm{\mF}^{l,1}\bm{p}^1\\
        \bm{\mF}^{l,2}\bm{p}^2\\
        \end{pmatrix},
    \end{equation}
    where each $\bm{\mF}^{l,i}, i\in\left\{1,2\right\}$ denotes a convolution with an interpolation kernel with small support.
    Denoting by $\bm{u}_j(\bm{\mF})$ the solution of the lower level problem \eqref{eq:bil_2} for the input image $\bm{f}_j$ in correspondence of a filter family $\bm{\mF}$, by $G(\bm{\A u},\bm{f}_j)=\frac{1}{2}\|\bm{\A u}-\bm{f}\|_2^2$ and by  $\mathcal{L}(\bm{\mF})=\frac{1}{sn}\sum_{j=1}^s \bm{\Phi}(\bm{u}_j(\bm{\mF}),\bm{g}_j)$, where $\bm{\Phi}(\bm{u}_j(\bm{\mF}),\bm{g}_j)=\frac{1}{2}\|\bm{u}_j(\bm{\mF})-\bm{g}_j\|_2^2$ the bilevel learning problem of finding the optimal interpolation filters $\bm{\mF}$ reads:
\begin{equation}\label{eq:bil_tv}
\begin{split}
&\min_{\bm{\mF}} ~\mathcal{L}(\bm{\mF}) + \mathcal{R}(\bm{\mF})\\
\text{s.t.}~\bm{u}_j\in\argmin_{\bm{u},\bm{q}}\max_{\bm{p}}\langle\bm{\mD}\bm{u}-&\bm{\mF}^*\bm{q}, \bm{p}\rangle + \lambda\|\bm{q}\|_{Z}+G(\bm{\A u},\bm{f}_j),\ j=1,\ldots,s
\end{split}
\end{equation}
where the regularization functional $\mathcal{R}$ is defined to impose prior constraints on the interpolation kernels. 
In particular, for $l=1,\ldots,L$, denoting by $\bm{\xi}^l=({\xi}_i^l)_i$ and $\bm{\eta}^l=({\eta}_i^l)_i$ the filter coefficients  of $\bm{\mF}^{l,1}$ and $\bm{\mF}^{l,2}$, respectively, we enforce that their sum has value $\mu\in\mathbb{R}$, or, as a shorthand notation that $\bm{\mF}\in (C_{\Sigma=\mu})^{L,2}$ where
$(C_{\Sigma=\mu})^{L,2}$ is the Cartesian product of $L$ copies of $(C_{\Sigma=\mu})^{1,2}$ defined as follows
\begin{equation}\label{eq:set_contraint}
    (C_{\Sigma=\mu})^{1,2}=\left\{\bm{\mF}^{l}:\ \sum_{i}\xi^l_i=\sum_{i}\eta^l_i=\mu, \right\}
\end{equation}
The functional $\mathcal{R}$ can thus be defined as:
\begin{equation}  \label{eq:regularizerF}
    \mathcal{R}(\bm{\mF})=\delta_{(C_{\Sigma=\mu})^{L,2}}(\bm{\mF})=\sum_{l=1}^L\delta_{(C_{\Sigma=\mu})^{1,2}}(\bm{\mF}^l),
\end{equation}
Solving problem \eqref{eq:bil_tv} may be very challenging due to the global non-convexity of the functional and its dependence on the solution of a non-smooth problem expressed in a primal-dual form.
A general method to approximately solve \eqref{eq:bil_tv} was proposed in \cite{chambolle2021learning}. It is reported in Algorithm \ref{alg:learn_tv}. Therein, the authors propose to use a proximal gradient method where the gradient of $\mathcal{L}(\bm{\mF})$ is computed by means of 
a linear approximation. More in details, since $\mathcal{L}$ is defined as the sum over the samples $j=1,\ldots,s$ in the dataset of the loss terms $\bm{\Phi}$, for a single sample $j$ the following approximation holds:
\begin{equation}\label{eq:gradient_approx}
    \nabla_{\bm{\mF}} ~\bm{\Phi}(\bm{u}_j(\bm{\mF}),\bm{g}_j)\approx -\left(\bm{Q}_j^K\otimes\bm{p}_j^K+\bm{q}_j^K\otimes\bm{P}_j^K\right)
\end{equation} 
where $\bm{q}_j^K,\bm{p}_j^K$ and $\bm{Q}_j^K, \bm{P}_j^K$ are, respectively, the last two terms of the saddle point and the corresponding adjoint states of the lower level problem  \eqref{eq:bil_tv} obtained by means of Algorithm \ref{alg:piggyback}, that is a piggyback primal dual algorithm \cite{griewank2003piggyback,Bogensperger2022} which jointly computes the solution both of the the lower level problem  \eqref{eq:bil_tv} and its associated biquadratic adjoint saddle-point problem.
The proximal operator of $\mathcal{R}(\bm{\mF})$, i.e. the projection over the set $(C_{\Sigma=\mu})^{L,2}$, can be derived following \cite{chambolle2021learning}. By separability \eqref{eq:regularizerF}, such projection can be computed separately for each filter projecting its weights onto the set $C_{\Sigma=\mu}$. For a given vector $\bm{\bar{x}}=(\bar{x}_1,\ldots,\bar{x}_n)\in\mathbb{R}^n$, the projection reads:
\begin{equation}\label{eq:mu_proj}
    \bm{\hat{x}}=\proj_{(C_{\Sigma=\mu})}(\bm{\bar{x}})\quad \Longleftrightarrow \quad \hat{x}_i=\bar{x}_i+\frac{\mu-\sum_{i=1}^n\bar{x}_i}{n},\quad i=1,\ldots n.
\end{equation}
The parameter $\mu$ is not set a priori but it is derived by a minimality argument as in \cite{chambolle2021learning}, solving the following problem:
\begin{equation}
    \min_{\mu,(\bm{x}^i)_{i=1}^m}\frac{1}{2}\sum_{i=1}^m\|\bm{x}^i-\bm{\bar{x}}^i\|^2,\quad\text{s.t. }\sum_{j=1}^n\bm{x}^i_j=\mu,\quad i=1,\ldots,m
\end{equation}
substituting \eqref{eq:mu_proj} in the previous formulation we get an explicit estimate of $\mu$ solving the following minimization problem:
\begin{equation}
    \hat{\mu}=\min_{\mu}\frac{1}{2}\sum_{i=1}^m\left(\mu-\sum_{j=1}^n \bm{\bar{x}}^i_j\right)^2\Longrightarrow \hat{\mu}=\frac{1}{m}\sum_{i=1}^m\sum_{j=1}^n \bm{\bar{x}}^i_j
\end{equation}

Note that imposing that the sum of all filter coefficients is equal to the same parameter $\mu$ allows to avoid the selection of an optimal regularization parameter $\lambda>0$, which depends both on the degradation and on the type of images considered. 
\begin{algorithm}[h]
\caption{Proximal gradient method to solve \eqref{eq:bil_tv}}\label{alg:learn_tv}
\begin{algorithmic}[0]
\State \textbf{Initialization}: choose $\bm{\mF}^{(0)}\in(C_{\Sigma=1})^{L,2}$, $\alpha>0$, $i=0$
\While {not converging}
    \State $\nabla_{\bm{\mF}} ~\bm{\Phi}(\bm{u}_j(\bm{\mF}^{(i)}),\bm{g}_j)\approx -\left(\bm{Q}_j^K\otimes\bm{p}_j^K+\bm{q}_j^K\otimes\bm{P}_j^K\right)$,\quad  $j=1,\ldots,s$
    \State $\nabla\mathcal{L}(\bm{\mF}^{(i)})=\frac{1}{sn}\sum_{j=1}^s \nabla_{\bm{\mF}} ~\bm{\Phi}(\bm{u}_j(\bm{\mF}^{(i)}),\bm{g}_j)$
    \State $\bm{\mF}^{(i+1)}=\proj_{(C_{\Sigma=\mu})^{L,2}}(\bm{\mF}^{(i)}-\alpha\nabla\mathcal{L}(\bm{\mF}^{(i)}))$
    \State $i=i+1$
\EndWhile
\State \textbf{Output}: learned interpolation kernel $\bm{\mF}^{(i+1)}$
\end{algorithmic}
\end{algorithm}

\begin{algorithm}[h]
\caption{Piggy-back primal-dual algorithm for solving low.lev. in \eqref{eq:bil_tv}.}\label{alg:piggyback}
\begin{algorithmic}[0]
\State \textbf{Initialization:}
$\bm{u}^0, \bm{U}^0\in\mathbb{R}^{n}$, 
$\bm{q}^0, \bm{Q}^0\in\mathbb{R}^{L\times 2n}$,
$\bm{p}^0, \bm{P}^0\in\mathbb{R}^{2n}$
\For{ $k=0,\ldots,K-1$ }
    \State 
    \begin{equation}
        \begin{cases}
        \bm{p}^{k+1}=\bm{p}^{k}+\sigma_{\bm{p}}(\bm{\mD}\bm{u}^k-\bm{\mF}^*q^k),&
        \bm{P}^{k+1}=\bm{P}^{k}+\sigma_{\bm{p}}(\bm{\mD}\bm{U}^k-\bm{\mF}^*\bm{Q}^k),\\
        \bar{\bm{p}}^{k+1}= \bm{p}^{k+1}+\theta(\bm{p}^{k+1}-\bm{p}^{k}), &
        \bar{\bm{P}}^{k+1}= \bm{P}^{k+1}+\theta(\bm{P}^{k+1}-\bm{P}^{k}), \\
        \bar{\bm{u}}^{k+1}= \bm{u}^{k}-\tau_{\bm{u}}\bm{\mD}^*\bar{\bm{p}}^{k+1}, &
        \bar{\bm{U}}^{k+1}= \bm{U}^{k}-\tau_{\bm{u}}(\bm{\mD}^*\bar{\bm{P}}^{k+1}+\nabla\bm{\Phi}(\bm{u}^k,g)),\\
        \bm{u}^{k+1} = \prox_{\tau_{\bm{u}} G \circ\bm{\A}}(\bar{\bm{u}}^{k+1}), &
        \bm{U}^{k+1} = \nabla\prox_{\tau_{\bm{u}} G \circ\bm{\A}}(\bar{\bm{u}}^{k+1})\cdot\bar{\bm{U}}^{k+1}, \\
        \bar{\bm{q}}^{k+1}= \bm{q}^{k}-\tau_{\bm{q}}\bm{\mF}\bar{\bm{p}}^{k+1}, &
        \bar{\bm{Q}}^{k+1}= \bm{Q}^{k}-\tau_{\bm{q}}\bm{F}\bar{\bm{P}}^{k+1},\\
        \bm{q}^{k+1}= \shrink_{\tau_{\bm{q}}\lambda}(\bar{\bm{q}}^{k+1}), & 
        \bm{Q}^{k+1}= \nabla\shrink_{\tau_{\bm{q}}\lambda}(\bar{\bm{q}}^{k+1})\cdot\bar{\bm{Q}}^{k+1}, \\
        \end{cases}
    \end{equation}
\EndFor
\State \textbf{Output}: approximate saddle point $(\bm{u}^K, \bm{q}^K, \bm{p}^K)$ and corresponding adjoint state $(\bm{U}^K, \bm{Q}^K, \bm{P}^K)$
\end{algorithmic}
\end{algorithm}

With the intent of testing \eqref{eq:bil_tv} for  different model operators $\bm{\A}$ (blur and super-resolution), we detail in the following the computation of the proximal operator $\prox_{\tau_{\bm{u}} G \circ\bm{A}}$ .

\textbf{Image Deblurring.} When $\bm{\A}=\mathbf{H}\in\mathbb{R}^{m\times n}$ is a convolution matrix, then
$\prox_{\tau_{\bm{u}} G \circ\mathbf{H}}$ can be computed by optimality:
\begin{equation}\label{eq:prox_debl}
    \hat{\bm{u}}=\mathrm{prox}_{\tau G \circ\bm{\mH}}(\bar{\bm{u}})\iff (\tau_{\bm{u}} \bm{\mH}^T\bm{\mH}+\bm{\mathrm{I}})\hat{\bm{u}}=\tau_{\bm{u}}\bm{\mH}^T\bm{g}+\bar{\bm{u}}
\end{equation}
The direct computation of $\hat{\bm{u}}$ thus requires the solution of an high-dimensional system. However, by assuming periodic boundary conditions, one can exploit the structure of the matrices involved (Block Circulant with Circulant Blocks) which allows a fast solution via the discrete Fast Fourier Transform (FFT) via the formula:
\begin{equation}
\hat{\bm{u}}=\mathcal{F}^{-1}\left(\frac{\tau_{\bm{u}}\overline{\mathcal{F}(\bm{\mH})}\mathcal{\bm{F}}(\bm{g})+\mathcal{\mF}(\bar{\bm{u}})}{\tau_{\bm{u}}\overline{\mathcal{F}(\bm{\mH})}\mathcal{F}(\bm{\mH})+1}\right)
\end{equation}
where $\mathcal{F}(\cdot)$ and $\overline{\mathcal{F}(\cdot)}$ denote the FFT and its conjugate, whereas $\mathcal{F}^{-1}(\cdot)$ is its inverse.

\textbf{Super-resolution.}
In the case of super-resolution, for $n=d^2 m, d>1$, the operator $\bm{\A}\in\R^{m\times n}$  maps high resolution images into their low resolution version. We consider here an operator $\bm{\A}=\mathbf{SH}$ where $\mathbf{H}\in\R^{n\times n}$ is a blur operator (as above) and $\mathbf{S}\in \left\{0,1\right\}^{m\times n}$ is the decimation operator which takes every $d$ columns/rows from the initial image to construct the new image. 
Its transpose $\bm{S}^T\in\mathbb{R}^{n\times m}$ is  the operator the interpolates the decimated images with zeros. 
For $\bm{\bar{u}}\in\R^n$, we have that by optimality, $\mathrm{prox}_{\tau_{\bm{u}} G \circ\bm{SH}}$ reads:
\begin{equation}\label{eq:prox_sr}
    \hat{\bm{u}}=\mathrm{prox}_{\tau_{\bm{u}} G \circ\bm{\mathrm{SH}}}(\bar{\bm{u}})\quad\iff\quad (\tau_{\bm{u}}\bm{\mH}^T\bm{\mathrm{S}}^T\bm{\mathrm{S}}\bm{\mH}+\bm{\mathrm{I}})\hat{\bm{u}}=\tau_{\bm{u}} \bm{\mH}^T\bm{\mathrm{S}}^T\bm{g}+\bar{\bm{u}}.
\end{equation}
Assuming period boundary conditions, a closed form to compute the solution of \eqref{eq:prox_sr} can be derived as in \cite{zhao2016fast} by 
exploiting the Woodbury formula and by factorizing the operator $\bm{\mathrm{S}}^T\bm{\mathrm{S}}$ as the Kronecker product of identity matrices and vector of ones \cite{Pragliola2023}. In the end, the solution reads as follows:
\begin{equation}
    \hat{\bm{u}}=\mathcal{F}^{-1}\left(\frac{\tau_{\bm{u}}\overline{\mathcal{F}(\bm{\mH})}\mathcal{\bm{F}}(\bm{\mathrm{S}}^T\bm{g})+\mathcal{F}(\bar{\bm{u}})}{\frac{\tau_{\bm{u}}}{d}\overline{\mathcal{F}(\bm{\mH})}\mathcal{F}(\bm{\mH})+1}\right).
\end{equation}
We have now all the ingredients to implement and carry out numerical simulations.

%% file: part_results.tex
In this section, we present and discuss the results obtained using  Algorithm~\ref{alg:learn_foe} and Algorithm \ref{alg:piggyback} presented in sections~\ref{sec:foe} and \ref{sec:TV_disc} for learning the FoE regularization model and the TV interpolation filters, respectively, on image deblurring and super-resolution problems.

\subsection{FoE bilevel learning for image deblurring}

We start evaluating the performance of Algorithm \ref{alg:learn_foe} estimating the parameters vector $\bm{\theta}=(\alpha_1,\ldots,\alpha_L, \bm{k_1},\ldots,\bm{k}_L)$ by solving the bilevel problem \eqref{eq:bil_foe}.  The routines used in this section are implemented in Python, using the PyTorch package for automatic differentiation routines. The experiment are performed on the Apocrita HPC server, using one GPU with eight cores with 7 GB of memory/core. 

\begin{figure}[!h]\centering
\begin{tabular}{@{}ccccccccc@{}}
\includegraphics[height=1.24cm]{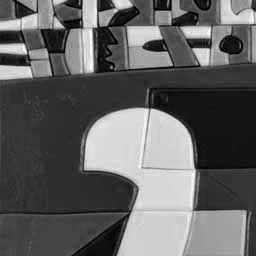}
&\includegraphics[height=1.24cm]{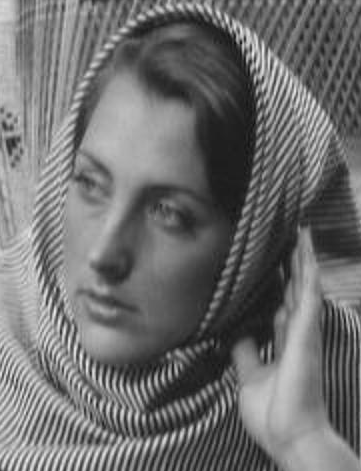}
&\includegraphics[height=1.24cm]{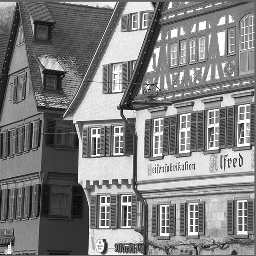}
&\includegraphics[height=1.24cm]{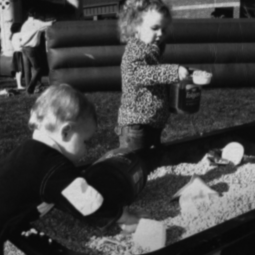}
&\includegraphics[height=1.24cm]{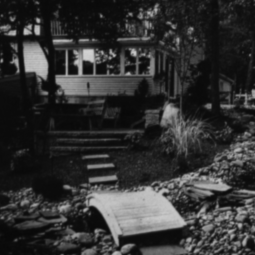}
&\includegraphics[height=1.24cm]{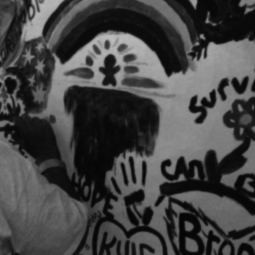}
&\includegraphics[height=1.24cm]{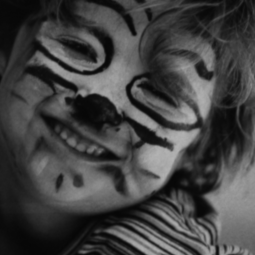}
&\includegraphics[height=1.24cm]{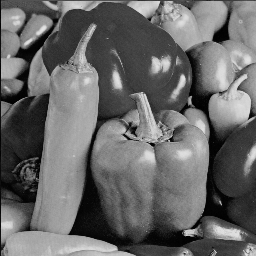}
&\includegraphics[height=1.24cm]{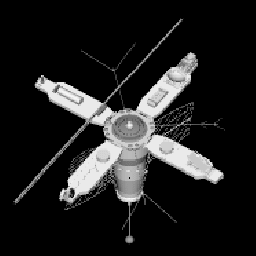}
\end{tabular}
\caption{Natural images used for training.}
\label{fig:train1}
\end{figure}

To generate the training set, we consider $N_{\text{TS}} = 9$ gray-scale natural images of different sizes (in particular, four images sized $255\times 255$, four sized $256\times 256$, and one sized $471\times 361$). These are shown in Figure \ref{fig:train1}.

\begin{figure}[h]
    \centering
    \begin{tabular}{@{}ccc@{\quad}ccc@{\quad}ccc@{}}
    \multicolumn{3}{c}{Gaussian blur}&
    \multicolumn{3}{c}{Disk blur}&
    \multicolumn{3}{c}{Motion blur}\\
        \includegraphics[height=1.24cm]{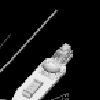} &\includegraphics[height=1.24cm]{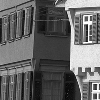} & \includegraphics[height=1.24cm]{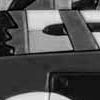} & \includegraphics[height=1.24cm]{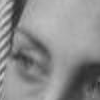} & \includegraphics[height=1.24cm]{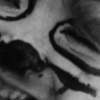} & \includegraphics[height=1.24cm]{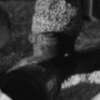} & \includegraphics[height=1.24cm]{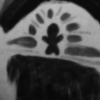} & \includegraphics[height=1.24cm]{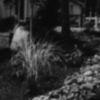} & \includegraphics[height=1.24cm]{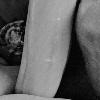}  \\
       \includegraphics[height=1.24cm]{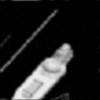} &\includegraphics[height=1.24cm]{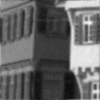} & \includegraphics[height=1.24cm]{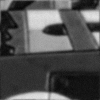} & \includegraphics[height=1.24cm]{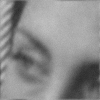} & \includegraphics[height=1.24cm]{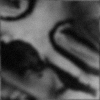} & \includegraphics[height=1.24cm]{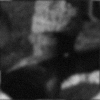} & \includegraphics[height=1.24cm]{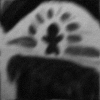} & \includegraphics[height=1.24cm]{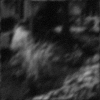} & \includegraphics[height=1.24cm]{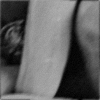}  
    \end{tabular}
    \caption{Examples of ground truth patches $\bm{g}_j$ (top) and degraded patches $\bm{f}_j$ (bottom)  included in the training set. Three different blur kernels are considered with additive Gaussian noise $\mathcal{N}(0,0.01^2\mathbf{I})$.}
    \label{fig:train_patch}
\end{figure}

From each test image,  $3$ patches of size $100\times 100$ are extracted.
The training set $\{(\bm{g}_j,\bm{f}_j)\}_{j=1}^s$
is thus composed of $s=9\times 3 = 27$ pairs of ground truth patches $\bm{g}_j$ and corrupted patches $\bm{f}_j$ obtained from $\bm{g}_j$ by means of three different blur kernels - namely a Gaussian blur with band 5, a disk blur with diameter 5 and a motion blur with length 5, 
all generated by assuming periodic boundary conditions
- and white Gaussian noise with standard deviation \texttt{sigma}$=0.01$.
Some exemplar ground truth and degraded patches
are shown in Figure \ref{fig:train_patch}.

\begin{figure}[b]
    \centering
    \begin{tabular}{@{}c@{\quad}c@{\quad}c@{}}
    \includegraphics[height=3cm]{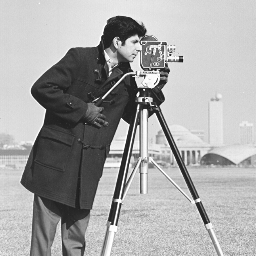}&   \includegraphics[height=3cm]{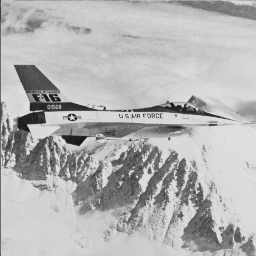}&   \includegraphics[height=3cm]{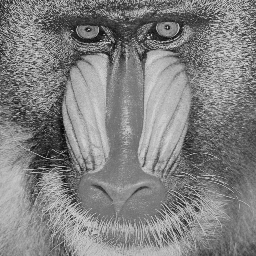}\\
    (a)&(b)&(c)
\end{tabular}
\caption{Test images: (a) \texttt{cameraman}, (b) \texttt{aircraft}, and (c) \texttt{mandrill}.}
\label{fig:or}
\end{figure}

\begin{figure}[tb]
    \centering
\begin{tabular}{@{}c@{\;\;}c@{\;\;}|@{\;\;}c@{\;}c@{\;}c@{\;}c@{}}
    \toprule 
    & Degraded 
    & \multicolumn{4}{c}{Reconstructions}\\[0.25em]
     \midrule
    & & $N_{\text{TS}} =9$ & $N_{\text{TS}} =30$ & $N_{\text{TS}} =50$ & $N_{\text{TS}} =75$ \\[0.25em]
       \rotatebox{90}{\quad Gaussian blur} &
       \includegraphics[height=2.25cm]{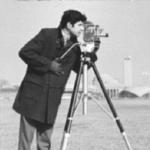} &
       \includegraphics[height=2.25cm]{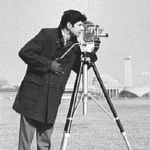} &
       \includegraphics[height=2.25cm]{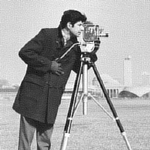} &
       \includegraphics[height=2.25cm]{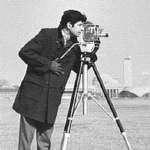} &
       \includegraphics[height=2.25cm]{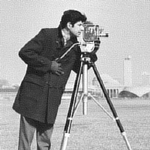} \\
       & PSNR = 24.48  & PSNR = 26.84  & PSNR = 26.94 & PSNR = 26.97 & PSNR = 26.97 \\[0.5em]
       \rotatebox{90}{\qquad Disc blur} &
       \includegraphics[height=2.25cm]{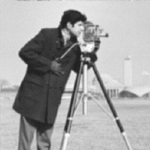} &
       \includegraphics[height=2.25cm]{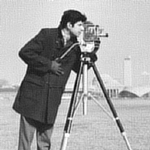} &
       \includegraphics[height=2.25cm]{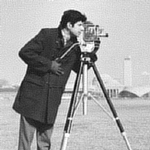} &
       \includegraphics[height=2.25cm]{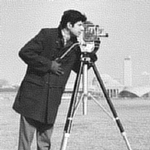} &
       \includegraphics[height=2.25cm]{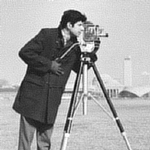}\\
       & PSNR = 23.28 & PSNR = 25.66  & PSNR = 25.79  & PSNR = 25.76 & PSNR = 25.79 \\[0.5em]
        \rotatebox{90}{\quad\; Motion blur} &
       \includegraphics[height=2.25cm]{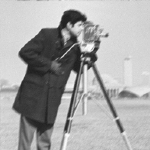} &
       \includegraphics[height=2.25cm]{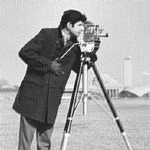} &
       \includegraphics[height=2.25cm]{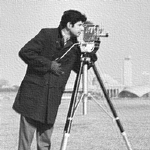} &
       \includegraphics[height=2.25cm]{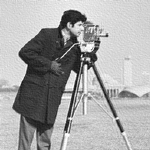} &
       \includegraphics[height=2.25cm]{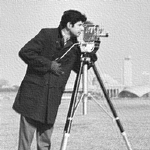} \\
       & PSNR = 23.60 & PSNR = 26.69  & PSNR = 26.76 & PSNR = 26.97 & PSNR = 26.82 \\
    \bottomrule
    \end{tabular}
    \caption{Reconstructed \texttt{cameraman} images for different types of blur (from top to bottom: Gaussian, disc, and motion blur), using different number of training samples (from left to right: $N_{\text{TS}} =9,\ 30,\ 50$, and $75$ images in the training dataset). Compare with ground truth in Figure~\ref{fig:or}(a).}
    \label{fig:test_patch_cameraman_comparison}
\end{figure}

\begin{figure}[!h]
    \centering
    \begin{tabular}{@{}c@{\;\;}c@{\;\;}|@{\;\;}c@{\;}c@{\;}c@{\;}c@{}}
    \toprule 
    & Degraded 
    & \multicolumn{4}{c}{Reconstructions}\\[0.25em]
    \midrule
    & & $N_{\text{TS}} =9$ & $N_{\text{TS}} =30$ & $N_{\text{TS}} =50$ & $N_{\text{TS}} =75$ \\[0.25em]
       \rotatebox{90}{\quad Gaussian blur} &
       \includegraphics[height=2.25cm]{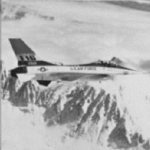} &
       \includegraphics[height=2.25cm]{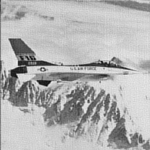} &
       \includegraphics[height=2.25cm]{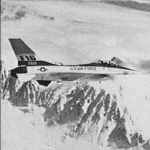} &
       \includegraphics[height=2.25cm]{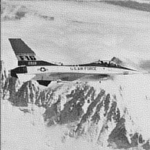} &
       \includegraphics[height=2.25cm]{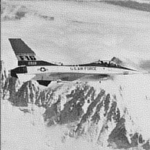} \\
       & PSNR = 26.41 & PSNR = 28.96 & PSNR = 29.04 & PSNR = 28.99 & PSNR = 29.08  \\[0.5em]
       \rotatebox{90}{\qquad Disc blur} &
       \includegraphics[height=2.25cm]{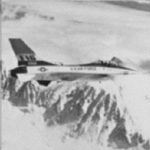} &
       \includegraphics[height=2.25cm]{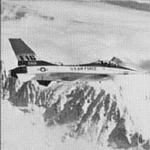} &
       \includegraphics[height=2.25cm]{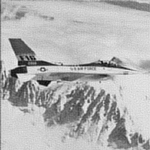} &
       \includegraphics[height=2.25cm]{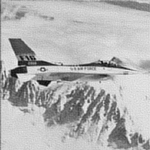} &
       \includegraphics[height=2.25cm]{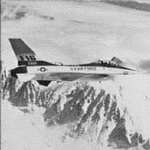}\\
       & PSNR = 25.21 & PSNR = 27.94 & PSNR = 27.94 & PSNR = 27.87 & PSNR = 27.95 \\[0.5em]
        \rotatebox{90}{\quad\; Motion blur} &
       \includegraphics[height=2.25cm]{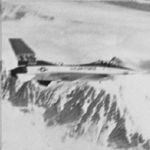} &
       \includegraphics[height=2.25cm]{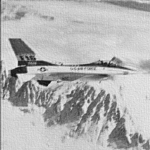} &
       \includegraphics[height=2.25cm]{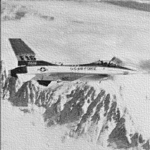} &
       \includegraphics[height=2.25cm]{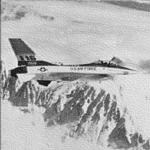} &
       \includegraphics[height=2.25cm]{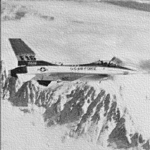} \\
       & PSNR = 25.09 & PSNR = 28.49 & PSNR = 28.45 & PSNR = 28.27 & PSNR = 28.47\\
    \bottomrule
    \end{tabular}
    \caption{Reconstructed \texttt{aircraft} images for different types of blur (from top to bottom: Gaussian, disc, and motion blur), using different number of training samples (from left to right: $N_{\text{TS}} = 9,\ 30,\ 50$, and $75$ images in the training dataset). Compare with ground truth in Figure~\ref{fig:or}(b).}
    \label{fig:test_patch_airplane_comparison}
\end{figure}

\begin{figure}[!h]
    \centering
    \begin{tabular}{@{}c@{\;\;}c@{\;\;}|@{\;\;}c@{\;}c@{\;}c@{\;}c@{}}
    \toprule 
    & Degraded 
    & \multicolumn{4}{c}{Reconstructions}\\[0.25em]
     \midrule
    & & $N_{\text{TS}} =9$ & $N_{\text{TS}} =30$ & $N_{\text{TS}} =50$ & $N_{\text{TS}} =75$ \\[0.25em]
       \rotatebox{90}{\quad Gaussian blur} &
       \includegraphics[height=2.25cm]{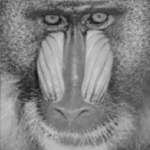} &
       \includegraphics[height=2.25cm]{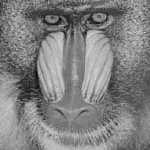} &
       \includegraphics[height=2.25cm]{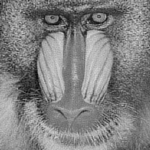} &
       \includegraphics[height=2.25cm]{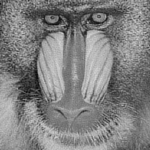} &
       \includegraphics[height=2.25cm]{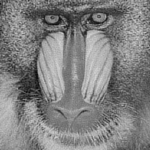} \\
       & PSNR = 21.67 & PSNR = 22.57 & PSNR = 22.76 & PSNR = 22.74 & PSNR =22.72 \\[0.5em]
       \rotatebox{90}{\qquad Disc blur} &
       \includegraphics[height=2.25cm]{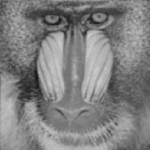} &
       \includegraphics[height=2.25cm]{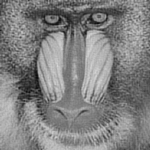} &
       \includegraphics[height=2.25cm]{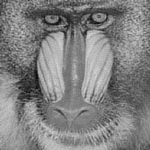} &
       \includegraphics[height=2.25cm]{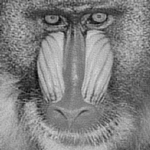} &
       \includegraphics[height=2.25cm]{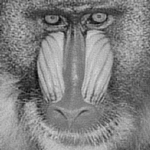}\\
       & PSNR = 20.76 & PSNR = 22.01 & PSNR = 22.22 & PSNR = 22.18 & PSNR = 22.17 \\[0.5em]
        \rotatebox{90}{\quad\; Motion blur} &
       \includegraphics[height=2.25cm]{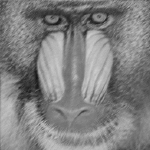} &
       \includegraphics[height=2.25cm]{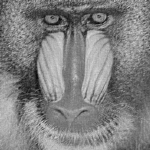} &
       \includegraphics[height=2.25cm]{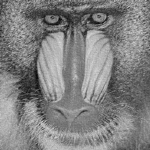} &
       \includegraphics[height=2.25cm]{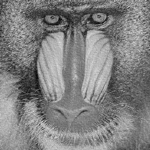} &
       \includegraphics[height=2.25cm]{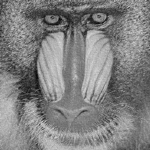} \\
       & PSNR = 21.24 & PSNR = 23.24 & PSNR = 23.60  & PSNR = 23.49 & PSNR = 23.54 \\
    \bottomrule
    \end{tabular}
    \caption{Reconstructed \texttt{mandrill} images for different types of blur (from top to bottom: Gaussian, disc, and motion blur), using different number of training samples (from left to right: $N_{\text{TS}} =9,\ 30,\ 50$, and $75$ images in the training dataset). Compare with ground truth in Figure~\ref{fig:or}(c).}
    \label{fig:test_patch_mandrill_comparison}
\end{figure}

\noindent The lower level problem in 
\eqref{eq:bil_foe} is solved via Algorithm \ref{alg:solve_lower}
with the following parameter setting: $\sigma = 10^{-4}$, $\beta = 0.5$, $\gamma_{\rm min} = 10^{-4}$, $\gamma_{\rm max}  = 1$. For the inner iterations a criterion based on relative error evaluated for the whole training set is used so that iterations are stopped when:
\begin{equation}\label{eq:rel_chg_inner}
  \frac{1}{s}\sum_{j=1}^s \frac{\|\bm{u}_{j,{(t+1)}}-\bm{u}_{j,{(t)}}\|_2}{\|\bm{u}_{j,{(t)}}\|_2}<\tau_{\mathrm{inner}} = 10^{-6}\,,\quad 
  \text{or}\quad t_{\max}=8000.
\end{equation}
The outer level problem is stopped 
at $k_{\max}=100$.
The learned filters and weights are tested on the restoration of three test 
images, namely \texttt{cameraman}, \texttt{aircraft} and \texttt{mandrill}, of size 256$\times$256 
shown in Figure \ref{fig:or}.
The test images considered are corrupted by same type of blur and additive white Gaussian noise as in the training phase.
The corrupted test images and the output restorations are shown in the two left-most columns of {Figures
\ref{fig:test_patch_cameraman_comparison},
\ref{fig:test_patch_airplane_comparison}, \ref{fig:test_patch_mandrill_comparison}}. We observe that in all restorations edges appear to be sharp and the noise is removed; nonetheless, in the case of motion blur the overall machinery tends to produce an over-deblurring effect, which is more evident in the case of \texttt{cameraman} and \texttt{aircraft} test images. 

To overcome such tendency, we performed the bilevel learning strategy described above with larger training datasets. In particular,  we used three subsets of the BSDS500 dataset \cite{MartinFTM01} consisting of $N_{\text{TS}} = 30$, $50$, and $75$ grey-scale images. Again, 3 patches of size $100 \times 100$ were extracted from each image and degraded as before. The results comparing different training set sizes are shown in the three right-most columns of Figures
\ref{fig:test_patch_cameraman_comparison},
\ref{fig:test_patch_airplane_comparison},  \ref{fig:test_patch_mandrill_comparison}, for the different test images. Moreover, in Table \ref{tab:foe_psnrAverage} we report the average PSNR values on both training and test data for the different blurs and sizes $N_{\text{TS}}$ of training samples. Notice that the average PSNR values for the training data are computed on patches of dimension $100\times 100$, while the average PSNR values on the testing data are computed on the whole $256\times 256$ images, thus explains the slight difference in the range of achieved values.

The output restorations and the average PSNR suggest that increasing the training set size does not significantly improve the deblurring results. On one hand, such behavior confirms the robustness of the bilevel set-up when dealing with the easiest scenarios of Gaussian and disk blur; on the other, the more challenging case of motion blur seems to require a slight modification of the lower-level problem adopted here. As an example, a global regularization parameter could be introduced, and added to the vector of learned parameters, so as to mitigate the sub-optimal scaling of the weights observed here.

\begin{figure}[h]
    \centering
    \includegraphics[height=4cm]{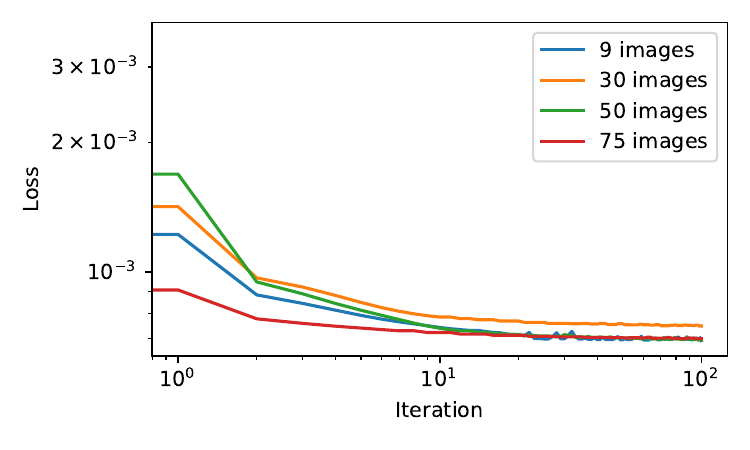}\\   
\caption{Log-log plot of the upper level loss function for different number $N_{\text{TS}}$ of images in the training sample.}
\label{fig:comparison-loss-different-datasets}
\end{figure}

Finally, Figure \ref{fig:comparison-loss-different-datasets} shows the log-log plot of the upper level loss function for the different training sets considered. We highlight that a faster decay of the loss is observed when larger datasets are used for training. Hence, there could be a trade off between the computational cost related to the processing of a training set and the number of iterations required for achieving a given tolerance on the loss function.

\begin{table}[t]
    \centering
    \begin{tabular}{@{}|l|l|c|c|c|c|c@{}}
    \cline{3-6}
    \multicolumn{1}{c}{} & & \; $N_{\text{TS}} =9$ \;  & \; $N_{\text{TS}} =30$ \; & \; $N_{\text{TS}} =50$ \; & \; $N_{\text{TS}} =75$ \; \\
    \hline
    \multirow{2}{*}{Gaussian blur} &  PSNR train & 31.15 & 29.75 & 29.52 & 29.37 \\
    \cline{2-6}
      &  PSNR test & 26.12 & 26.25 & 26.23 & 26.26 \\
    \hline
    \multirow{2}{*}{Disc blur} &  PSNR train & 30.98 & 29.44 & 29.79 & 30.20 \\
    \cline{2-6}
      &  PSNR test & 25.21 & 25.31 & 25.27 & 25.31 \\
    \hline
    \multirow{2}{*}{Motion blur} &  PSNR train & 30.41 & 29.22 & 29.82 & 29.70 \\
    \cline{2-6}
      &  PSNR test & 26.14 & 26.27 & 26.15 & 26.28 \\
      \hline
    \end{tabular}
    \caption{Average PSNR on both training and test data for the different deblurring problems (from top to bottom: Gaussian, disc and motion blur) and number of training samples (from left to right: $N_{\text{TS}} =9,\ 30,\ 50$, and $75$ images in the training dataset). }
    \label{tab:foe_psnrAverage}
\end{table}


\subsection{TV discretization learning for image deblurring and super-resolution}

In this section, we evaluate the performance of the bilevel learning strategy developed in Section \ref{sec:TV_disc} for solving problem \eqref{eq:bil_tv} computing optimal discretization filters $\bm{\mF}$ for TV regularization  via Algorithm \ref{alg:piggyback}. Several numerical results are reported in which both the training and test datasets are let vary with different degradation settings. The routines used in this section are implemented in Python, using the PyTorch package to exploit automatic differentiation routines. The experiment are performed on a Dell PowerEdge server, equipped with a Nvidia Tesla V100 having 32 GB of memory. 

\begin{figure}[h]
\centering
\begin{adjustbox}{width=\columnwidth,center}
\begin{tabular}{ccccccccc}
\frame{\includegraphics[height=1.24cm]{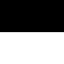}}&\frame{\includegraphics[height=1.24cm]{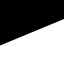}}&\frame{\includegraphics[height=1.24cm]{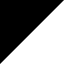}}&\frame{\includegraphics[height=1.24cm]{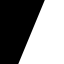}}&\frame{\includegraphics[height=1.24cm]{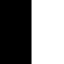}}&\frame{\includegraphics[height=1.24cm]{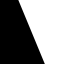}}&\frame{\includegraphics[height=1.24cm]{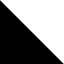}}&\frame{\includegraphics[height=1.24cm]{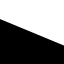}}\\
\frame{\includegraphics[height=1.24cm]{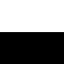}}&\frame{\includegraphics[height=1.24cm]{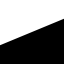}}&\frame{\includegraphics[height=1.24cm]{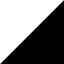}}&\frame{\includegraphics[height=1.24cm]{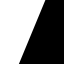}}&\frame{\includegraphics[height=1.24cm]{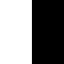}}&\frame{\includegraphics[height=1.24cm]{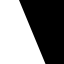}}&\frame{\includegraphics[height=1.24cm]{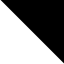}}&\frame{\includegraphics[height=1.24cm]{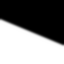}}
\end{tabular}
\end{adjustbox}
\caption{Examples of training images considered  for learning the optimal discretization of TV.}
\label{fig:train2}
\end{figure}

For all experiments we consider a small $2\times 2$ pixels neighborhood, whose associated convolution kernels $\bm{\mF}^{l,i}$, $l=1,\ldots, L$, $i\in\left\{1,2\right\}$ have size $2\times 3$ and $3\times 2$ for the horizontal and vertical components of the dual variable $\bm{p}$, respectively. A different numbers of filters weights $L\in\{2,3,4,8\}$ is considered, so as to analyze how the reconstruction quality varies depending on $L$. 
To model possible symmetries present in the observed data, we allow the possibility to incorporate some invariances in the filter weights. In particular, we allow transpose symmetry for the cases $L=2,3$, and rotational symmetry of angles $\frac{\pi}{2}$ for the cases $L=4,8$. The projection onto these symmetries groups can be found in \cite{chambolle2021learning}. The training data consists of $s=64$ images of size $64\times 64$ pixels with ground truth being nothing but binary images with an edge with equi-spaced orientations $\theta_{j}=2j\pi/s$, $j=0,\ldots,s-1$. A subset of these images is shown in Figure \ref{fig:train2}. To enrich the dataset and simulate partial volume effects, the dataset includes a small random shift of the discontinuities from the center. The test dataset is generated in a similar way. 

We use $K=2000$ iterations of Algorithm \ref{alg:piggyback} to compute the derivatives of the loss function at each outer iteration. In Algorithm \ref{alg:learn_tv} we use a constant stepsize $\alpha=100$  and run the learning algorithm for $I=500$ iterations, since we empirically observe that the objective function stabilizes its decrease around that number. 
We used a warm-start strategy for the adjoint states variables to improve the accuracy for gradients approximation of the lower-level problem.

In the following, we compare the results of the proposed strategy for two different imaging tasks (namely, deblurring and super-resolution) with those obtained using two handcrafted discretization filters. In particular, we consider the horizontal and vertical forward differences, which we denote by FD, and the filters proposed by Condat in \cite{condat2017discrete}, which we denote by CD3 and CD4 where the numbers represent how many directional filters are considered. 
Results are compared in terms of their PSNR averaged over the considered dataset.

\subsubsection{Image deblurring}
We consider the problem of learning optimal TV discretization using as training images the ones in Figure \ref{fig:train2}. In order to avoid unwanted artifacts on the images boundaries image padding  using reflexive condition is used and a crop is performed to get their initial size. Three datasets are considered with Gaussian blur of varied width quantified by the standard deviation of the Gaussian kernel $\varsigma$:
\begin{itemize}
    \item GaussianA: $\varsigma=0.5$ (small blur).
    \item GaussianB: $\varsigma=1$ (medium blur).
    \item GaussianC: $\varsigma=1.5$ (high blur).
\end{itemize}
The results of the method are summarized in Table \ref{tab:0_noise_blur} both for the training and the test set. Since PSNR is only slightly lower on the test data, we deduce that overfitting is quite limited. By imposing symmetries onto the filter weights leads to better results outperforming their non-symmetric counterpart, as it reads in particular for the GaussianA setting. On the other hand, using a larger number of filters does not seem to improve significantly the quality of the results.
In Table \ref{tab:filter_noisefree_deblur} we report the learned filters for the GaussianA setting with no noise. Color coding has to be interpreted as follows: black corresponds to the lowest filter weight, white to the maximum filter weight and the gray levels in between correspond to intermediate values between the two extremes with various intensities. The filters $L=2$ and $L=3$ show a mild symmetry with respect to the horizontal axis. 
 
\begin{table}[h]
    \centering
    \begin{adjustbox}{width=\columnwidth,center}
    \begin{tabular}{c|c|c|c}
    \hline
    $L=2$ & $L=2$ (s) & $L=3$ & $L=3$ (s) \\
    \hline
         \includegraphics[width=0.2\textwidth]{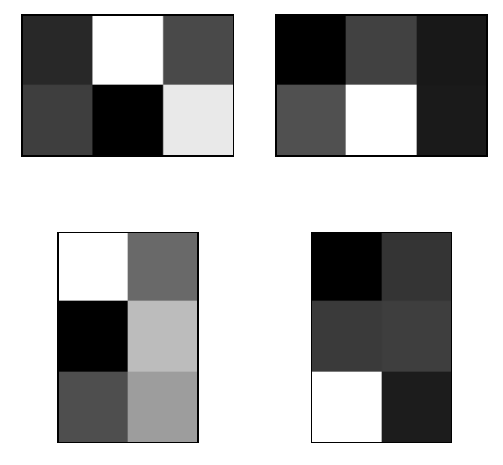} &
         \includegraphics[width=0.2\textwidth]{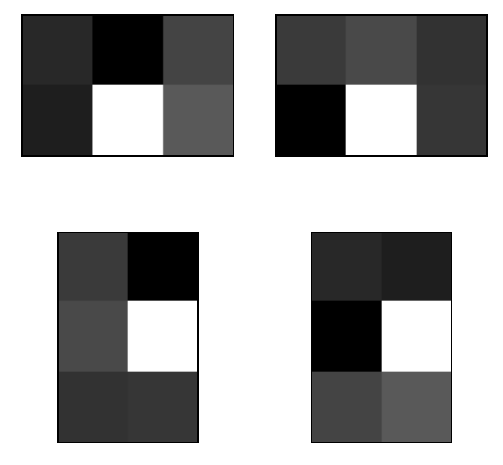} &
         \includegraphics[width=0.3\textwidth]{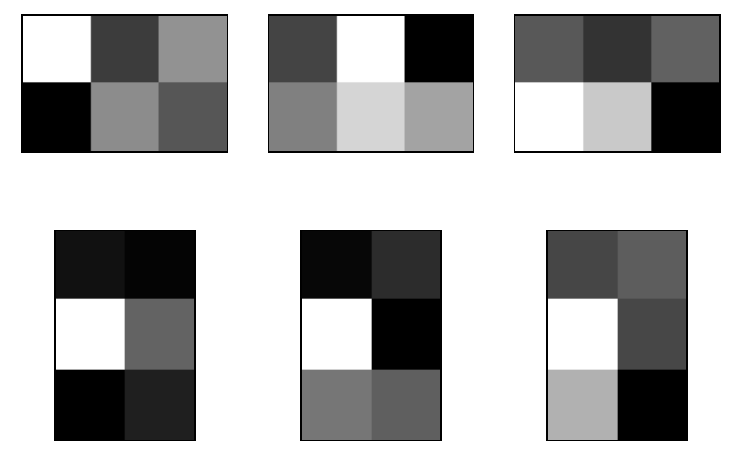} &
         \includegraphics[width=0.3\textwidth]{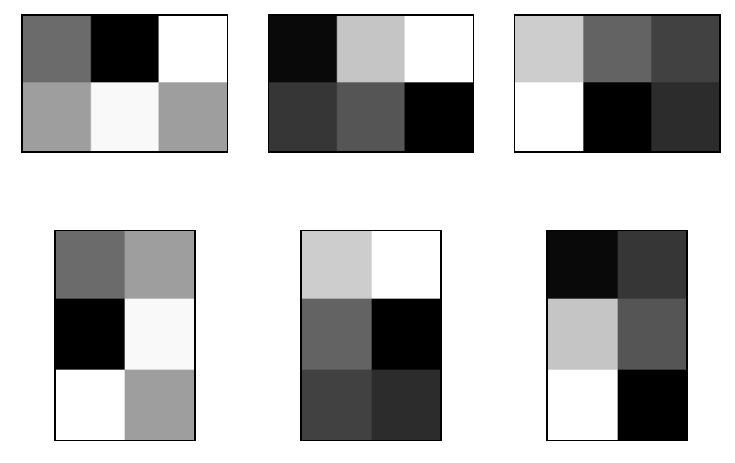} 
         \\
    \end{tabular}
    \end{adjustbox}
    \centering
    \begin{adjustbox}{width=\columnwidth,center}
    \begin{tabular}{c|c}
    \hline
    $L=4$ (s) & $L=8$ (s) \\
    \hline
         \includegraphics[width=0.33\textwidth]{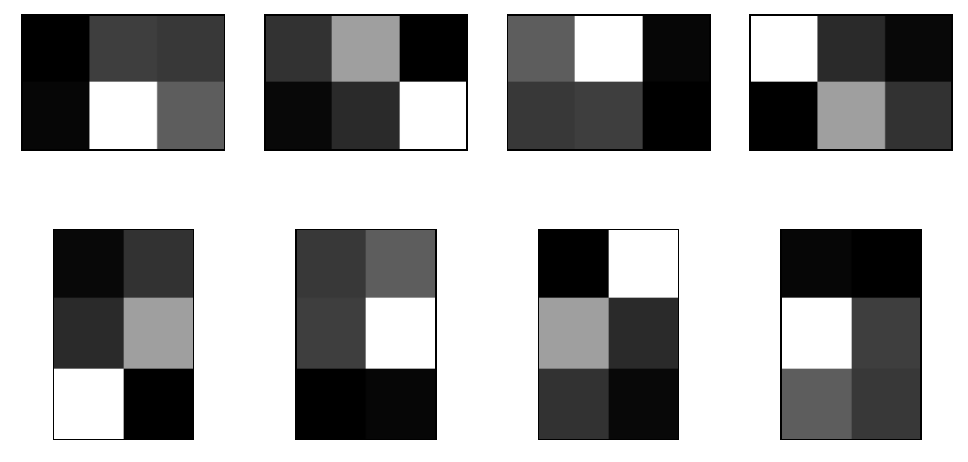} &
         \includegraphics[width=0.67\textwidth]{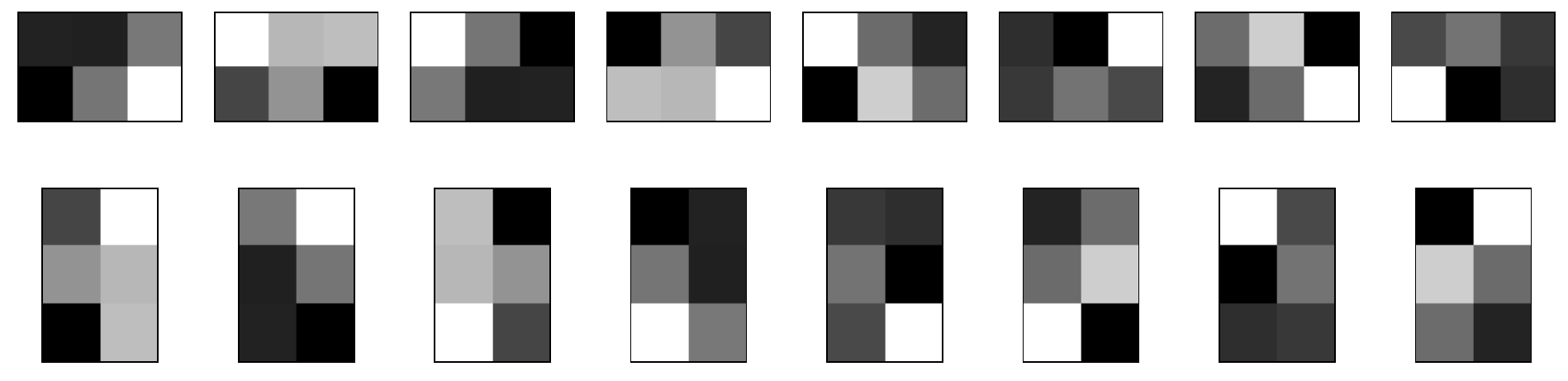} 
         \\ \hline
    \end{tabular}
    \end{adjustbox}
    \caption{Learned filters for the noise-free GaussianA ($\varsigma=0.5$) setting. Symmetric filters are indicated by (s).}
    \label{tab:filter_noisefree_deblur}
\end{table}
\begin{table}[h]
    \centering
    \begin{adjustbox}{width=\columnwidth,center}
    \begin{tabular}{@{}|l|l|ccc|cccccc|@{}}
\cline{3-11}
\multicolumn{1}{}{} & &
 FD & CD3 & CD4 & $L=2$ & $L=2$ (s) & $L=3$ & $L=3$ (s) & $L=4$ (s) & $L=8$ (s) \\
\hline

\multirow{2}{*}{GaussianA}& PSNR train &    39.12 &      43.1 &     43.78 &                 45.30 &                    45.55 &                45.58 &                    44.27 &                    45.33 &                    44.03 \\
\cline{2-11} 
& PSNR test  &    39.16 &     43.09 &     43.78 &                44.93 &                    45.29 &                 45.30 &                    44.06 &                    44.95 &                    43.28 \\
\hline
\multirow{2}{*}{GaussianB}& PSNR train &     33.60 &     39.76 &     40.08 &                38.52 &                    41.53 &                39.64 &                    39.27 &                    41.28 &                    41.51 \\
\cline{2-11} 
& PSNR test  &    33.63 &     39.62 &     40.01 &                38.42 &                    41.15 &                39.16 &                     39.30 &                    40.86 &                     40.9 \\
\hline
\multirow{2}{*}{GaussianC} & PSNR train &    31.38 &     37.51 &     37.56 &                 38.50 &                    38.89 &                38.73 &                    37.99 &                    37.14 &                    39.19 \\
\cline{2-11} 
& PSNR test  &    31.35 &     37.55 &     37.58 &                38.58 &                    38.85 &                38.42 &                    38.09 &                    37.24 &                    39.02 \\
\hline

\end{tabular}
\end{adjustbox}
    \caption{Comparison between average PSNR of hand-crafted and learned filters evaluated on both training and test data for the different noise-free deblurring problems.}
    \label{tab:0_noise_blur}
\end{table}

We then compared the results with the ones obtained for the same degradation settings GaussianA, GaussianB and GaussianC with additional white Gaussian noise with standard deviation \texttt{sigma}$=0.01$. A summary of the result can be found in Table \ref{tab:1e-2_noise_blur}. Similar  considerations to those same discussed in the noise-free case can be drawn.

\begin{table}[!h]
    \centering
    \begin{adjustbox}{width=\columnwidth,center}
    \begin{tabular}{@{}|l|l|ccc|cccccc|@{}}
\cline{3-11}
\multicolumn{1}{}{} & & 
 FD & CD3 & CD4  & $L=2$ & $L=2$ (s) & $L=3$ & $L=3$ (s) & $L=4$ (s) & $L=8$ (s) \\
\hline

\multirow{2}{*}{GaussianA}& PSNR train &  39.07 &     42.91 &     43.56 &                41.06 &                    45.08 &                43.88 &                    45.23 &                    44.49 &                    45.02 \\
\cline{2-11}
& PSNR test  &    39.05 &     42.91 &     43.58 &                 41.00 &                    44.87 &                43.36 &                    44.93 &                    43.93 &                    44.43\\
\hline
\multirow{2}{*}{GaussianB}& PSNR train &    33.61 &     39.52 &      39.90 &                39.32 &                    38.84 &                 40.00 &                     40.90 &                     40.40 &                    41.08 \\
\cline{2-11}
& PSNR test  &    33.57 &     39.53 &     39.88 &                39.18 &                     38.5 &                39.72 &                    40.73 &                    40.36 &                     40.80\\
\hline
\multirow{2}{*}{GaussianC} & PSNR train &    31.38 &     37.46 &     37.52 &                37.82 &                    37.95 &                37.81 &                    38.46 &                    37.76 &                    38.13 \\
\cline{2-11}
& PSNR test  &    31.37 &     37.47 &     37.54 &                37.84 &                    38.02 &                37.84 &                    38.31 &                    37.84 &                    38.12 \\
\hline
\end{tabular}
\end{adjustbox}
    \caption{Comparison between average PSNR of hand-crafted and learned filters evaluated on both training and test data for the different deblurring + noisy problems. }
    \label{tab:1e-2_noise_blur}
\end{table}

\begin{figure}[h]
    \centering
    \includegraphics[width=\textwidth]{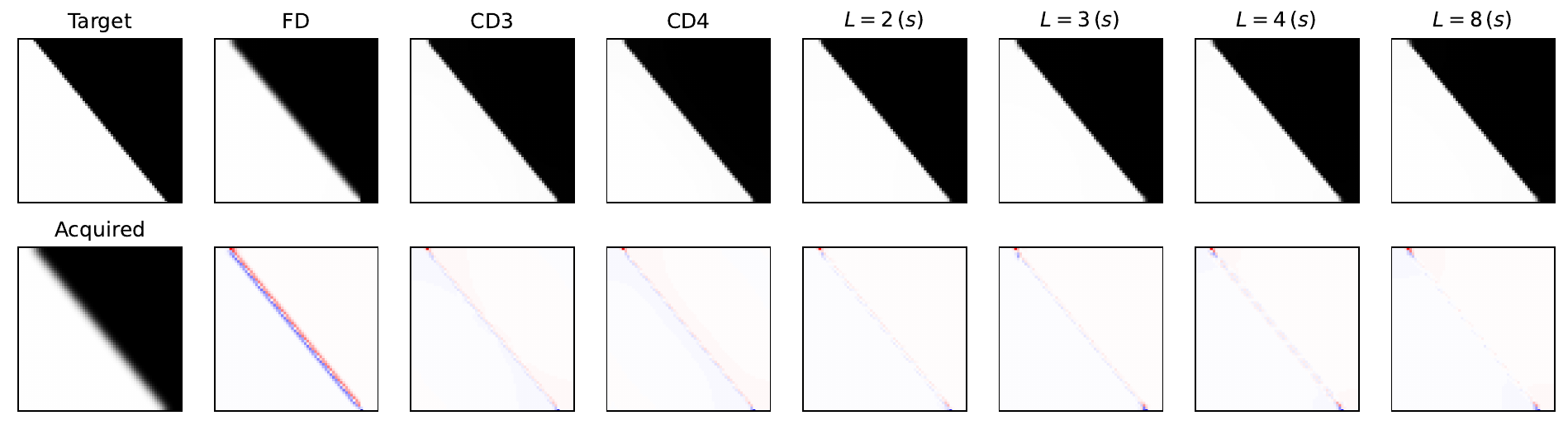}
    \caption{Results (top row) and color-coded error map (bottom row) highlighting differences for ground truth image corresponding to $\theta\approx7\pi/4$ in noise-free GaussianC setting.}
    \label{fig:0_noise_blur_examples}
\end{figure}

\begin{figure}[h]
    \centering
    \includegraphics[width=\textwidth]{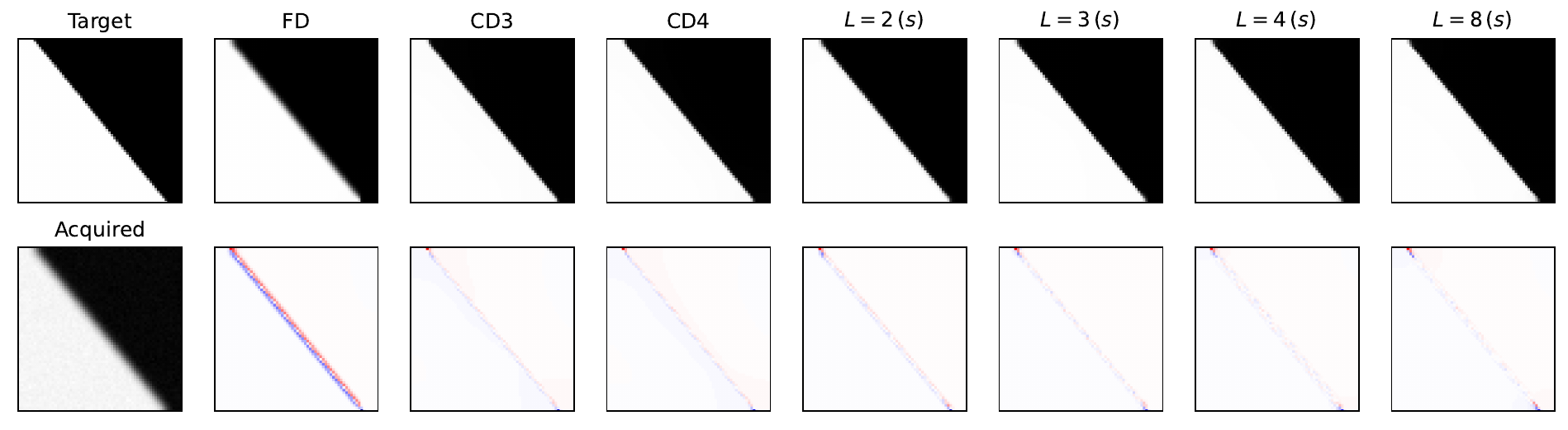}
    \caption{Results (top row) and color-coded error map (bottom row) highlighting differences for ground truth image corresponding to $\theta\approx7\pi/4$ in noisy GaussianC setting.}
    \label{fig:1e-2_noise_blur_examples}
\end{figure}

For the most challenging blur setting (GaussianC), we report in Figure \ref{fig:0_noise_blur_examples}-\ref{fig:1e-2_noise_blur_examples} false-color error plots comparing the reconstruction obtained with the corresponding  target, both for the noise-free and noisy case, respectively. Red represents here the pixel error value 1 and blue -1, while white depicts an pixel error near zero. It is evident that the handcrafted FD filters  shows significant biases along edges. On the contrary, the CD filters perform almost as well as the learned filters, which give the best results. In particular, the setting $L=8$ (s), in Figure \ref{fig:0_noise_blur_examples}, restores the discontinuity with nearly zero error.


\subsubsection{Image super-resolution}

We now consider the problem of super-resolution (SR) for the same dataset of images described above. To do so, we consider as a forward operator a column- and row-wise downsampling operator of factor $d=2$, which amounts to discarding every second column/row. To assess the quality of the learned filters w.r.t.~the resolution loss only, we consider in the modeling a composition with a convolution matrix $\mathbf{H}$ corresponding to a small point spread function ($\varsigma=0.1$) coinciding, upon discretization, with a Dirac delta.
The bilevel learning procedure is repeated for these data.
The results for both the training and the test set are summarized in Table \ref{tab:noblur_sf2}. Compared to the deblurring case, the difference in terms of PSNR between the training and the test set is more evident, meaning that, to some extent, we are overfitting the training set. For this problem, the incorporating symmetries on the filter weights does not provide any competitive improvements on the image quality.
In Table \ref{tab:filter_noisefree_SR} we report the learned filters on the noise-free SR setting. Differently form the deblurring problem, the filters with $L=2$ and $L=3$ do not show any natural symmetry. 

\begin{table}[]
    \centering
    \begin{adjustbox}{width=\columnwidth,center}
    \begin{tabular}{c|c|c|c}
    \hline
    $L=2$ & $L=2$ (s) & $L=3$ & $L=3$ (s) \\
    \hline
         \includegraphics[width=0.2\textwidth]{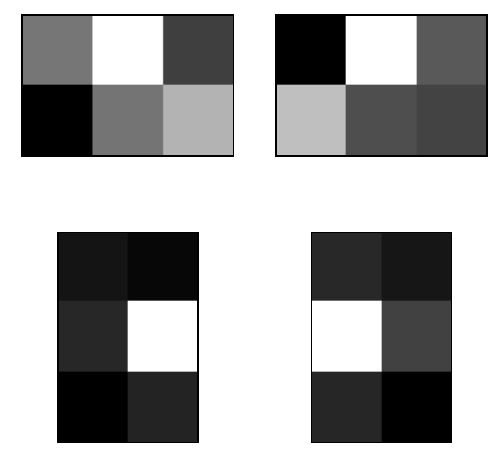} &
         \includegraphics[width=0.2\textwidth]{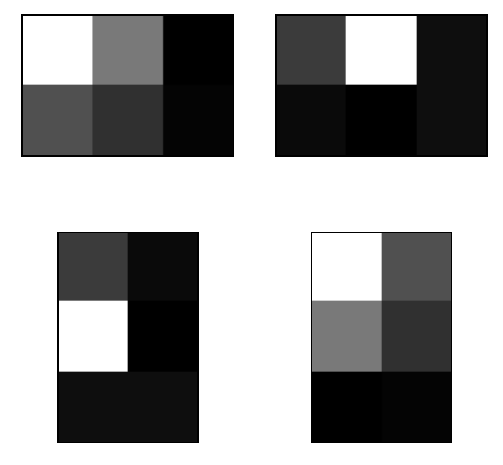} &
         \includegraphics[width=0.3\textwidth]{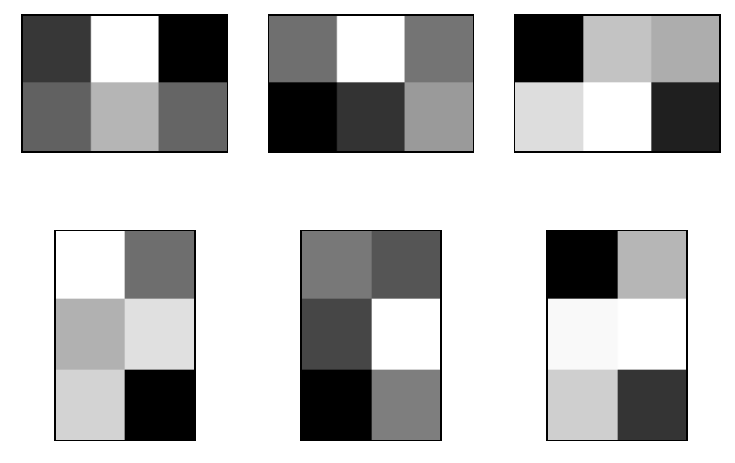} &
         \includegraphics[width=0.3\textwidth]{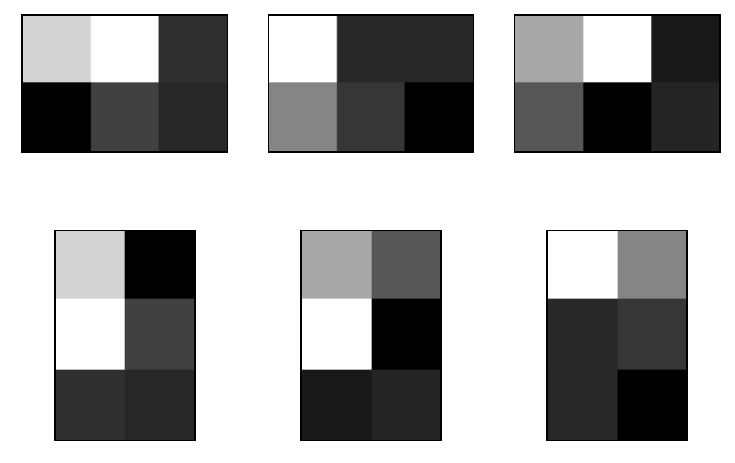} 
         \\
    \end{tabular}
    \end{adjustbox}
    \centering
    \begin{adjustbox}{width=\columnwidth,center}
    \begin{tabular}{c|c}
    \hline
    $L=4$ (s) & $L=8$ (s) \\
    \hline
         \includegraphics[width=0.33\textwidth]{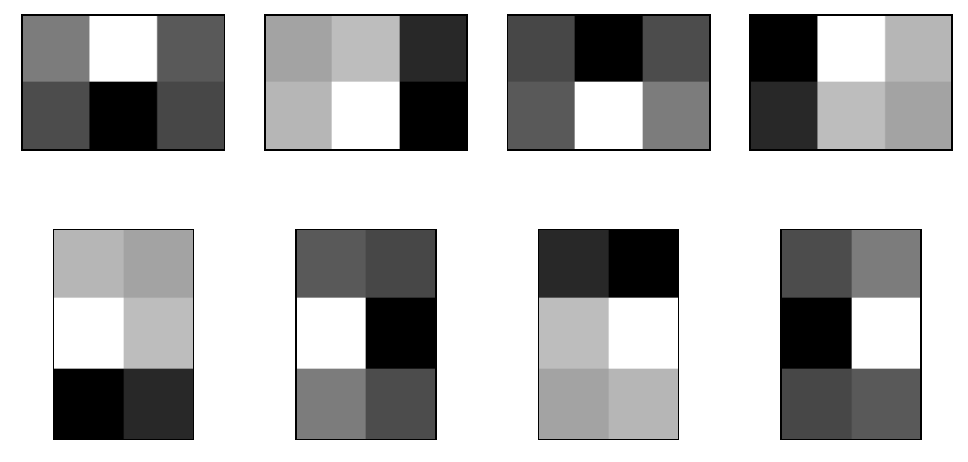} &
         \includegraphics[width=0.67\textwidth]{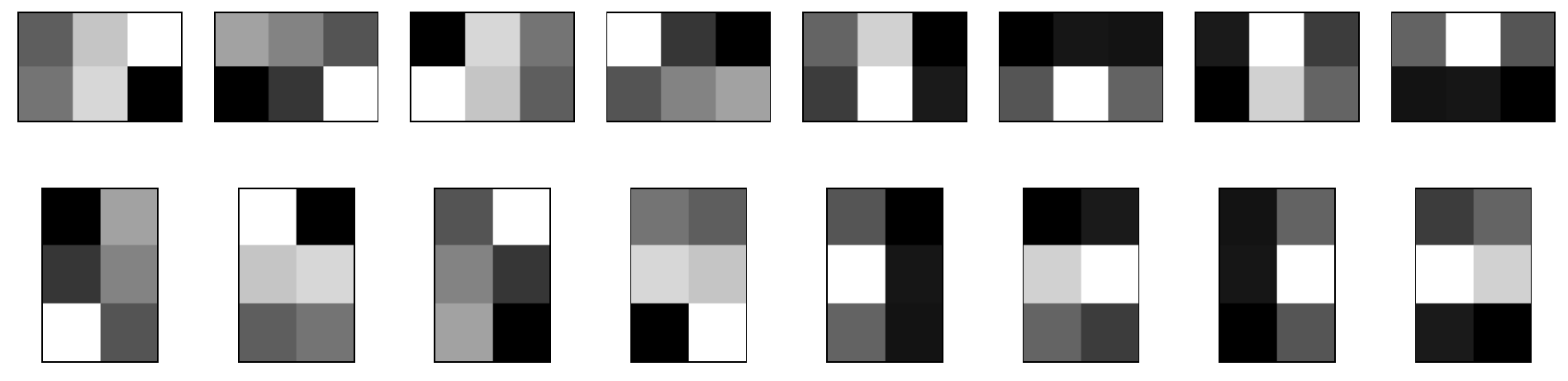} 
         \\
         \hline
    \end{tabular}
    \end{adjustbox}
    \caption{Learned filters for the super-resolution noise-free setting. Symmetric filters are indicated by (s).}
    \label{tab:filter_noisefree_SR}
\end{table}

\begin{table}[h]
    \centering
    \begin{adjustbox}{width=\columnwidth,center}
    \begin{tabular}{|l|l|ccc|cccccc|}
\cline{3-11}
\multicolumn{1}{}{} & & FD & CD3 & CD4 & $L=2$ & $L=2$ (s) & $L=3$ & $L=3$ (s) & $L=4$ (s) & $L=8$ (s) \\
\hline
\multirow{2}{*}{Noisefree} &
PSNR train &    28.11 &     30.47 &     30.04 &                 33.80 &                    33.83 &                33.85 &                    32.95 &                    32.76 &                    33.03 \\
\cline{2-11} 
& PSNR test  &    27.91 &     30.07 &     29.65 &                31.92 &                     32.00 &                32.63 &                    31.79 &                    32.72 &                    32.35 \\
\hline
\multirow{2}{*}{Noisy} &
PSNR train &    27.76 &      30.20 &     29.76 &                32.93 &                    31.96 &                 33.10 &                    33.32 &                    33.14 &                    32.33 \\
\cline{2-11} 
& PSNR test  &    28.19 &      30.40 &     30.04 &                 32.60 &                    31.57 &                32.49 &                    33.26 &                    32.46 &                    31.98 \\
\hline
\end{tabular}
\end{adjustbox}
    \caption{PSNR of handcrafted and learned filters evaluated on both the training and test data for the noise-free SR setting with $d=2$.}
    \label{tab:noblur_sf2}
\end{table}

\begin{figure}
    \centering
    \includegraphics[width=\textwidth]{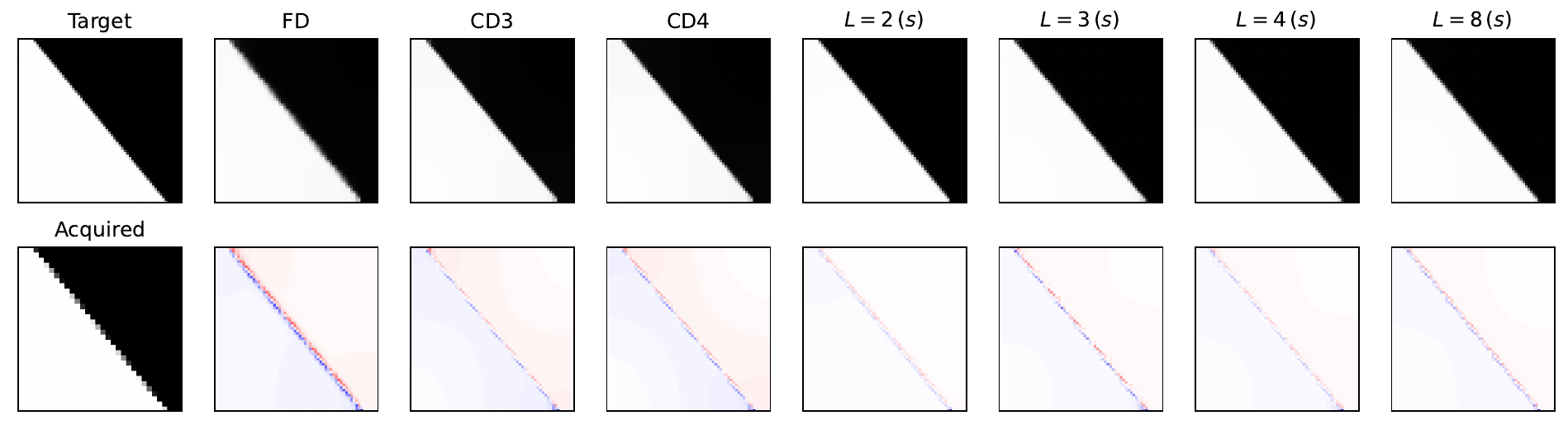}
    \caption{Results (top row) and color-coded error map (bottom row) highlighting differences for ground truth image corresponding to $\theta\approx7\pi/4$ in  noise-free SR setting.}
    \label{fig:0_noise_sr_examples}
\end{figure}

\begin{figure}
    \centering
    \includegraphics[width=\textwidth]{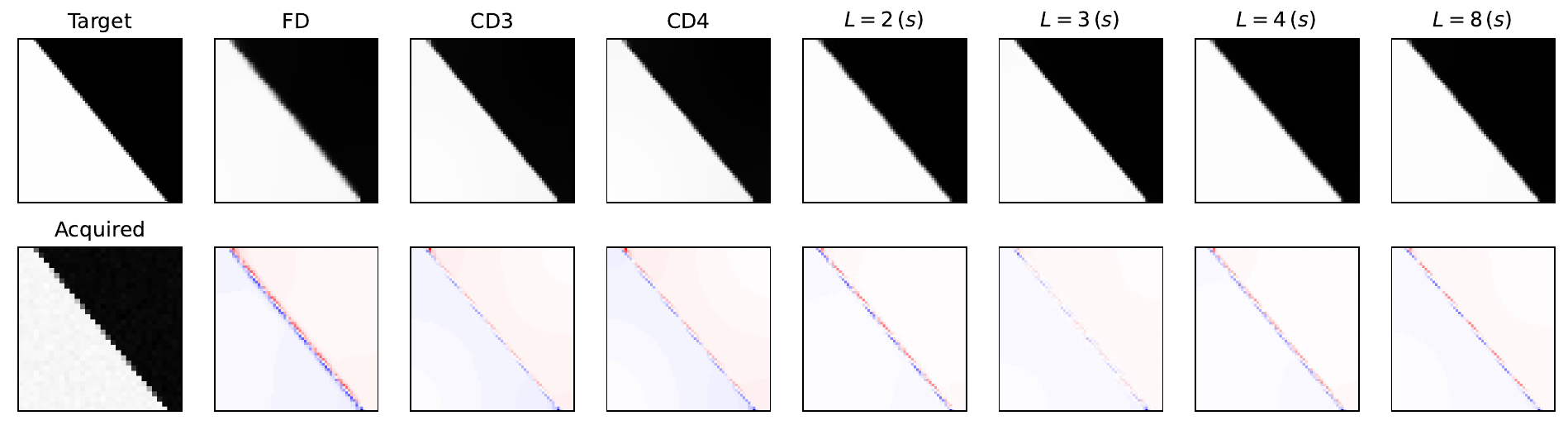}
    \caption{
    Results (top row) and color-coded error map (bottom row) highlighting differences for ground truth image corresponding to $\theta\approx7\pi/4$ in  noisy SR setting.}
    \label{fig:1e-2_noise_sr_examples}
\end{figure}

Similarly to deblurring case, we report in Figure \ref{fig:0_noise_sr_examples}-\ref{fig:1e-2_noise_sr_examples} false-color error plots between the reconstructed image with learned filters with respect to the corresponding target image. The handcrafted FD filters  provide the worst results, 
with several discretization biases along the edge. Regarding the learned filters, the $L=2$ (s) setting visually provides the best restoration for both the noise-free and noisy case, where the $L=4$ (s) configuration gives a good result as well.


\subsubsection{Crossover Testing}

We conclude our numerical exploration by investigating whether the filters learned on a specific task are able to generalize to a different task. In particular, we focus on four different settings and consider whether filters learned for a specific noise-free/noisy deblurring (GaussianB) and SR problem, generalize well to the other tasks.
We further compare the filters with $L=2$ (s) and $L=3$ (s), which were the ones delivering overall good results across the experiments above. In Table \ref{tab:crossover_L2} and \ref{tab:crossover_L3} we report the confusion matrices (w.r.t. PSNR) of the results obtained in comparison with the handcrafted filters (CD4). It is evident that the filters obtained for the SR task generalize very well on the deblurring task. In general, for $L=2$ (s) learning the filters on a task without noise seems to give filters that generalize on the problem also in presence of the noise. For the case $L=3$ (s) it happens the opposite, namely the filters learned on a noisy task generalize well on the corresponding noiseless task.

\begin{table}[h]
    \centering
    \begin{adjustbox}{width=\columnwidth,center}
    \begin{tabular}{c|l|cc|cc|c|}
        \multicolumn{2}{c}{} & \multicolumn{4}{c}{Learning task} & \multicolumn{1}{c}{Handcrafted}\\
        \cline{3-7}
        \multicolumn{2}{c|}{}  & GaussianB & GaussianB noisy & SR & SR noisy & CD4\\
        \cline{2-7}
        \multirow{4}{*}{Evaluation Task\quad } & GaussianB & 41.15 & 38.60 & 38.55 & 35.49 & 40.01 \\ \cline{2-7}
         & GaussianB noisy & 40.82 &    38.50 & 38.20 &   35.28 & 39.88\\ \cline{2-7}
         & SR & 28.93 &    30.73 &  32.00 &    31.52 & 29.65\\ \cline{2-7}
         & SR noisy & 29.34 &    30.81 & 32.47 &    31.57 & 30.04 \\
        \cline{2-7}
    \end{tabular}
	\end{adjustbox}
    \caption{Crossover testing results for $L=2$ (s).}
    \label{tab:crossover_L2}
\end{table}

\begin{table}[h]
    \centering
    \begin{adjustbox}{width=\columnwidth,center}
    \begin{tabular}{c|l|cc|cc|c|}
        \multicolumn{2}{c}{} & \multicolumn{4}{c}{Learning task} & \multicolumn{1}{c}{Handcrafted}\\
        \cline{3-7}
        \multicolumn{2}{c|}{}  & GaussianB & GaussianB noisy & SR & SR noisy & CD4\\
        \cline{2-7}
        \multirow{4}{*}{Evaluation Task\quad } & GaussianB & 39.30 &    41.02 & 37.16 &    37.87 & 40.01 \\ \cline{2-7}
         & GaussianB noisy &   38.82 &   40.73 & 36.90 &   37.79 & 39.88\\ \cline{2-7}
         & SR & 28.26 &   30.14 & 31.79 &   32.58  & 32.58\\ \cline{2-7}
         & SR noisy & 28.49 &    30.57 & 31.91 &    33.26 & 30.04 \\
        \cline{2-7}
    \end{tabular}
	\end{adjustbox}
    \caption{Crossover testing results for $L=3$ (s).}
    \label{tab:crossover_L3}
\end{table}
As previously observed in \cite{chambolle2021learning}, we observe that the learned filters highly depend on the specific inverse problem used for the learning phase.